\documentclass{amsart}

\usepackage{amsmath}
\usepackage{amscd}

\theoremstyle{plain}
\newtheorem{thm}{Theorem}[section]
\newtheorem{lem}[thm]{Lemma}
\newtheorem{prop}[thm]{Proposition}
\newtheorem{reduc}[thm]{Reduction}

\theoremstyle{definition}

\theoremstyle{remark}
\newtheorem{rem}{Remark}

\numberwithin{equation}{section}

\newcommand{\rl}{cabled linkage}
\newcommand{\ul}{{\mathcal L}}
\newcommand{\fdf}[1]{{\emph{ #1}}}
\newcommand{\conf}{{\mathcal C}}
\newcommand{\sconf}{{\mathcal S}{\mathcal C}}
\newcommand{\vrtcs}{{\mathcal V}}
\newcommand{\edge}{{\mathcal E}}
\newcommand{\qas}{quasialgebraic set}

\newcommand{\TS}{Tarski-Seidenberg}
\newcommand{\mand}{\text{ and }}
\newcommand{\mif}{\text{ if } }
\newcommand{\mfor}{\text{ for } }

\newcommand{\qf}{quasifunctional}
\newcommand{\qfl}{quasifunctional linkage}
\newcommand{\rd}{restricted domain}
\newcommand{\compose}{\circ}

\newcommand{\rh}[2]{\rho _{\ul #1,#2}}
\newcommand{\capshun}[2]{\caption{ #2}\label{fig:#1}}
\newcommand{\commentout}[1]{}
\newcommand{\figsix}{section \ref{sec:4.5}}

\newcommand{\Eu}{\mathrm {Euc}}

\newcommand{\Tran}{\mathrm {Tran}}
\newcommand{\lineseg}[2]{\overline{#1#2}}

       \input boxedeps.tex 
       \input boxedeps.cfg 
\newcommand{\Myepsf}[1]{{\BoxedEPSF {#1}}}

\begin{document}

\title{Configuration Spaces of Linkages in ${\mathbb  R}^n$}
\author[KING]{Henry C. King}
\address{
Department of Mathematics\\
University of Maryland\\
College Park, Maryland, 20742}
\email{hking@math.umd.edu}
\subjclass{Primary 14P05, 14P10; Secondary 57R99,14P20}

\date{November 23, 1998.}

\begin{abstract}
This paper studies the configuration space of  all possible positions of a linkage in 
${\mathbb R}^n$.
For example, it shows that for every compact algebraic set,
there is a linkage whose configuration space is 
analytically isomorphic to a finite number of copies
of the algebraic set.
If flexible edges are allowed, any compact set given by polynomial
equalities and inequalities is the configuration space of a linkage.
This paper also studies semiconfiguration spaces 
of all possible positions of a finite number of points on a linkage.
For example any compact semialgebraic set is such a semiconfiguration space.
\end{abstract}

\maketitle

\section{Linkages}

Loosely speaking, a linkage is an ideal mechanical device
consisting of a bunch of stiff rods sometimes attached at their ends
by rotating joints. A realization of a linkage in ${\mathbb  R}^n$ is some way of
placing this linkage in ${\mathbb  R}^n$.
The configuration space for a linkage is the space of all such realizations,
which can be determined by looking at all possible positions
of the ends of all the rods.
A semiconfiguration space of a linkage is the space of all possible positions
 of only some of ends of the rods, we ignore the other ends.
  For example, what figure does a particular point on the linkage
  trace out?

In this paper we will give characterizations of configuration spaces
and semiconfiguration spaces of linkages as well as of \rl s,
for $n\ge 3$.
(In a \rl\ you also can attach flexible cables between rods.)
The characterizations of these spaces for $n=2$,
planar linkages, was studied in \cite{KM}, \cite{K1}, and \cite{K2}.
The results for $n\ge 3$ turn out to be analogous to the $n=2$
results, but in a couple of places the proofs are different.
In particular, we can completely characterize
semiconfiguration spaces, we can characterize 
configuration spaces of \rl s up to analytic isomorphism,
and we can characterize configuration spaces of linkages
up to analytically trivial finite covers.

Let us now define linkages more precisely.
Suppose $L$ is a finite one dimensional simplicial
complex, in other words, a finite set $\vrtcs(L)$
of vertices and a finite set $\edge(L)$ of edges
between certain pairs of vertices.
An \fdf{abstract linkage} is a finite one dimensional simplicial
complex $L$ with a mapping $\ell \colon \edge(L)\to (0,\infty )$.
You should think of $\ell $ as giving the length of each edge.
A \fdf{realization} of an abstract linkage $(L,\ell )$ in ${\mathbb  R}^n$ is a mapping
$\varphi \colon \vrtcs(L)\to {\mathbb  R}^n$ so that 
$|\varphi (v)-\varphi (w)|=\ell (\lineseg{v}{w})$ 
for all  edges $\lineseg {v}{w}\in \edge(L)$.

We will often wish to fix some of the vertices of a linkage whenever we
take a  realization.
So we say that a \fdf{classical linkage in ${\mathbb  R}^n$}  is a foursome $\ul=(L,\ell ,V,\mu )$ where
 $(L,\ell )$ is an abstract linkage, $V\subset \vrtcs(L)$ is a subset of its vertices, and $\mu \colon V\to {\mathbb  R}^n$.
 So $V$ is the set of fixed vertices and $\mu $ tells where to fix them.
The configuration space of realizations is defined by:
$$\conf(\ul)=\left\{\,\varphi \colon \vrtcs(L)\to {\mathbb  R}^n\Biggm|
\begin{array}{cll}
\varphi (v)=\mu (v)&\mif\ v\in V\\
   |\varphi (v)-\varphi (w)|
   =\ell (\lineseg{v}{w})&\mfor \text{all edges } \lineseg {v}{w}\in \edge(\ul)
   \end{array}
   \,\right\}$$

A \fdf{\rl\ in ${\mathbb  R}^n$}\ is a quintuple $(L,\ell ,V,\mu ,F)$ where
$(L,\ell ,V,\mu )$ is a classical linkage and $F\subset \edge(L)$.
We will think of the edges in $F$ as being flexible
rather than rigid.
A physical model for such a \rl\ would
 have the edges in $\edge(L)-F$ be rigid rods 
 but the edges in $F$ are just ropes or cables
connecting two vertices.
Thus in a realization, two vertices connected by an edge $e$ in
$F$ would only be constrained to have distance $\le \ell (e)$.
The configuration space of a \rl\ is given by:
$$\conf(\ul)=\left\{ \varphi \colon \vrtcs(L)\to {\mathbb  R}^n \Biggm|
\begin{array}{cll}
\varphi (v)=\mu (v) &\mif v\in V,\\
  |\varphi (v)-\varphi (w)|\le \ell (\lineseg{v}{w})&\mif\ \lineseg {v}{w}\in F  \\
   |\varphi (v)-\varphi (w)|=\ell (\lineseg{v}{w})&\mif\ \lineseg v{w}\in \edge(L)-F
\end{array}
\right\}$$

From now on, the word linkage will refer to a \rl\ in ${\mathbb  R}^n$.
If we wish to refer to a linkage without any flexible edges,
we will call it a classical linkage.
If $F$ is empty we get a classical linkage.

If $W\subset \vrtcs(\ul)$ is a collection of vertices
of a linkage $\ul$, the \fdf{semiconfiguration space} is
the set of restrictions to $W$ of realizations of $\ul$,
$$\sconf(\ul,W)=\{\,\varphi \colon W\to {\mathbb  R}^n 
\text{ so that } \varphi =\varphi '|_W
\text{ for some } \varphi '\in \conf(\ul)  \,\}$$
Thus $\sconf(\ul,W)$ keeps track of the positions of only those vertices in $W$
and ignores the positions of other vertices.

To characterize configuration and semiconfiguration spaces,
we need a few definitions.
A \fdf{real algebraic set} is the set of solutions of a collection of
real polynomial equations in some ${\mathbb  R}^m$.
We define a \fdf{\qas}\ to be a subset of ${\mathbb  R}^m$ of the form
 $$\{\,x\in {\mathbb  R}^m\,\mid\, p_i(x)=0,\, i=1,\ldots ,\ell \mand\ q_j(x)\ge 0,\, j=1,\ldots ,k\,\}$$
 for some polynomials $p_i$ and $q_j$.
Finally, a semialgebraic set is a finite union of differences of \qas s.
In other words, a semialgebraic set is a finite union of sets of the form
$$\{\,x\in {\mathbb  R}^m\mid p_i(x)=0, q_j(x)\ge 0, \mand r_k(x)>0\,\}$$
for collections of polynomials $p_i$, $q_j$ and $r_k$.
Real algebraic sets and semialgebraic sets are well studied,
 but I am not aware of any literature on \qas s.

We will use two notions of isomorphism.
If $X\subset {\mathbb  R}^k$ and $Y\subset {\mathbb  R}^m$ then we say a 
homeomorphism $f\colon X\to Y$ is an \fdf{isomorphism}
if $f$ and $f^{-1}$ are both restrictions of entire rational functions,
for example polynomials.
We say $f$ is an \fdf{analytic isomorphism}
if $f$ and $f^{-1}$ are both restrictions of analytic maps,
i.e., maps locally given by power series.
So any isomorphism is analytic, but the converse is not true.
All analytic isomorphisms in this paper will actually be 
polynomials in one direction, but the inverse may involve
 square roots of positive quantities, and hence only be analytic.

\begin{rem}
Note that  any configuration
space $\conf(\ul)$ of a classical linkage
is an algebraic set in $({\mathbb  R}^n)^{\vrtcs(\ul)}$,
since it is the solutions of the polynomial equations
$|y_i-y_j|^2=\ell_{ij}^2$ and $y_i=z_i$ for $i\in V$.
Likewise, the configuration space of a \rl\ is a
\qas\ in $({\mathbb  R}^n)^{\vrtcs(\ul)}$.
Finally, a semiconfiguration space $\sconf(\ul,W)$
of a linkage $\ul$ is a semialgebraic set in $({\mathbb  R}^n)^W$.
To see this, note that it is the image of the \qas\
$\conf(\ul)$ under the projection map 
$({\mathbb  R}^n)^{\vrtcs(\ul)}\to ({\mathbb  R}^n)^W$.
But by the \TS\ theorem \cite{S},
the projection of a semialgebraic set is semialgebraic.
\end{rem}

We will prove the following theorems characterizing configuration spaces of
classical and \rl s, for $n\ge 3$.
For $n=2$, these theorems were proven in \cite{KM} and \cite{K1}.

\begin{thm}\label{thm:1}
Suppose $\ul$ is a classical linkage in ${\mathbb  R}^n$, $n\ge 2$.
Then $\conf(\ul)$
is isomorphic to $X\times ({\mathbb  R}^n)^k$ for some compact real algebraic set $X$.
The integer $k$ is the number of connected components of $\ul$
with no fixed vertices.

Conversely, if $X$ is a compact real algebraic set and $k\ge 0$,
there is a classical linkage $\ul$ and a finite set $F$ so that
$\conf(\ul)$ is analytically isomorphic to $X\times F\times ({\mathbb  R}^n)^k$.
In fact, there is  an analytic function
$\alpha \colon X\times F\times ({\mathbb  R}^n)^k\to ({\mathbb  R}^n)^m$ so that 
$$\conf(\ul)=\{\, (x,\alpha (x,f,y)) 
\mid x\in X , y\in ({\mathbb  R}^n)^k, \mand\ f\in F\,\}$$
and so that the map $(x,y,f)\mapsto (x,\alpha (x,f,y))$ is an analytic isomorphism.
\end{thm}

\begin{thm}\label{thm:2}
Suppose $\ul$ is a \rl\ in ${\mathbb  R}^n$, $n\ge 2$.
Then $\conf(\ul)$
is isomorphic to $X\times ({\mathbb  R}^n)^k$ for some compact \qas\ $X$.
The integer $k$ is the number of connected components of $\ul$
with no fixed vertices.

Conversely, if $X$ is a compact \qas\ and $k\ge 0$,
there is a \rl\ $\ul$ so that
$\conf(\ul)$ is analytically isomorphic to $X\times ({\mathbb  R}^n)^k$.
In fact, there is  an analytic function
$\alpha \colon X\times ({\mathbb  R}^n)^k\to ({\mathbb  R}^n)^m$ so that 
$$\conf(\ul)=\{\, (x,\alpha (x,y)) 
\mid x\in X, \mand\ y\in ({\mathbb  R}^n)^k\,\}$$
and so that the map $(x,y)\mapsto (x,\alpha (x,y))$ is an analytic isomorphism.
\end{thm}

We need a few more definitions before stating the characterization
of semiconfiguration spaces, since that characterization is
more precise.

Let $\Eu(n)$ denote the group of Euclidean motions of ${\mathbb  R}^n$.
So a general element of $\Eu(n)$ is of the form
$z\mapsto Q(z)+z_0$  where
$Q \in O(n)$ is an orthogonal matrix and $z_0\in {\mathbb  R}^n$ is a constant.
We say a subset $Z\subset ({\mathbb  R}^n)^k$ is \fdf{virtually compact}
if either $Z$ is compact, or $Z$ is invariant under the diagonal action
of $\Eu(n)$, with compact quotient.


\begin{thm}\label{thm:3}
Suppose $X\subset ({\mathbb  R}^n)^k$, $n\ge 2$.
Then the following are equivalent:
\begin{enumerate}
\item There is a \rl\  $\ul$ in ${\mathbb  R}^n$ and a $W\subset \vrtcs(\ul)$
so that $\sconf(\ul,W)=X$.
\item There is a classical linkage $\ul$ in ${\mathbb  R}^n$ and a $W\subset \vrtcs(\ul)$
so that $\sconf(\ul,W)=X$.
\item After perhaps permuting the coordinates,
 $X=Y_1\times  Y_2 \times \cdots \times Y_m $
where each 
$Y_i\subset ({\mathbb  R}^n)^{k_i}$ is  a virtually compact semialgebraic set.
\end{enumerate}
\end{thm}

Perhaps it is useful to restrict attention to the compact case.
Then any compact semialgebraic set is the semiconfiguration space
of a classical linkage.
Any compact \qas\ is analytically isomorphic to the configuration 
space of a \rl.
For any compact real algebraic set $X$ there is a classical linkage whose
configuration space is analytically isomorphic to a number of disjoint
copies of $X$.

As a byproduct of the proof of Theorem \ref{thm:1},
the cardinality of $F$ will be $2^b$.
This comes about because there are $b$ vertices which each have
two distinct positioning modes.
Essentially, Theorem \ref{thm:2} is proven by tethering these $b$
vertices to fixed points so that only one of the two modes is
possible (although other vertices are tethered as well).

\begin{rem}
If $\ul=(L,\ell ,V,\mu ,F)$ is a linkage and $\beta \in \Eu(n)$ we may form a
linkage $\beta (\ul)=(L,\ell ,V,\beta \compose \mu ,F)$ by applying $\beta $ to the image all
fixed vertices.  Since $\beta $ preserves distances, we know that
$\beta (\conf(\ul))=\conf(\beta (\ul))$ and $\beta (\sconf(\ul,W))=\sconf(\beta (\ul),W)$.
This is what we refer to as translating and rotating the linkage.
We may also rescale $\ul$ as follows.
If $\lambda $ is a positive number, then $\lambda \ul=(L,\lambda \ell ,V,\lambda \mu ,F)$.
We have $\lambda \conf(\ul)=\conf(\lambda \ul)$ and $\lambda \sconf(\ul,W)=\sconf(\lambda \ul,W)$.
\end{rem}

Because of the above remark, (semi)configuration spaces of
linkages with few fixed points have a great deal of symmetry.
We get more precise characterization theorems by taking this into account.
While complete for semiconfiguration spaces,
this characterization is incomplete for configuration spaces.
For simplicity, we only state these results for connected linkages.
Using Lemma \ref{lem:6.2} below, one could then formulate
 analogous results for nonconnected linkages.

To fix notation,
let $e_i\in {\mathbb  R}^n$ denote the unit vector whose only nonzero component
is a 1 in the $i$-th place.
For $1\le k\le n$ we have a subgroup $O(k)\subset \Eu(n)$ which we fix on as
the set of $\beta \in \Eu(n)$ so that $\beta (0)=0$ and $\beta (e_i)=e_i$ for all $i\le n-k$.
We let $\Tran(n)\subset \Eu(n)$ denote the subgroup of translations,
maps of the form $z\mapsto z+z_0$.


	\begin{thm}\label{thm:4x}
	Suppose $\ul$ is a connected linkage with exactly $m$ fixed vertices
and $W\subset \vrtcs(\ul)$,
then:
\begin{enumerate}
\item  If $m=0$ then $\sconf(\ul,W)$ and $\conf(\ul)$ are invariant under
the action of $\Eu(n)$, with compact quotient.
\item  If $1\le m\le n$, then $\sconf(\ul,W)$ and $\conf(\ul)$ are compact 
and invariant 
under a subgroup of $\Eu(n)$ conjugate to $O(n-m+1)$.
\item  To make part 2 above sharper, suppose $m\ge 1$ and
$T$ is an affine subspace of
${\mathbb  R}^n$ which contains the images of all fixed vertices.
Then $\sconf(\ul,W)$ and $\conf(\ul)$ are compact and invariant 
under the subgroup of elements $\Eu(n)$ which fix 
all points of $T$.
(This subgroup is conjugate to $O(n-\dim T)$.)
\end{enumerate}
	\end{thm}

	\begin{thm}\label{thm:4}
Let $Z\subset ({\mathbb  R}^n)^k$ be a virtually compact semialgebraic set, 
$n\ge 2$,
and suppose that $Z$ is invariant under the diagonal action of a
subgroup $G$ of $\Eu(n)$, with compact quotient.
Then there is a connected classical linkage $\ul$ and a $W\subset \vrtcs(\ul)$ so that
$\sconf(\ul,W)=Z$ and so that
\begin{enumerate}
\item If $G=\Eu(n)$,
then  $\ul$ has no fixed vertex.
\item  If $G$ is conjugate to $O(m)$, $1\le m\le n$,
then  $\ul$ has $n-m+1$ fixed vertices,
and these vertices are fixed at points on the fixed subspace of $G$.
\item  Otherwise,  $\ul$ has only $n+1$ fixed vertices.
\end{enumerate}
\end{thm}

	\begin{thm}\label{thm:5}
Let $Z\subset ({\mathbb  R}^n)^k$ be a compact algebraic set, $n\ge 2$,
and suppose that $Z$ is invariant under the diagonal action of a
subgroup $G$ of $\Eu(n)$.
Then there is a connected classical linkage $\ul$ and a finite set $F$ so that
$\conf(\ul)$ is analytically isomorphic to $Z\times F$ and so that
\begin{enumerate}
\item   If $G$ is conjugate to $O(1)$,
then  $\ul$ has $n$ fixed vertices,
and these vertices are fixed at points on the fixed subspace of $G$.
\item   Suppose $G$ is conjugate to $O(2)$
and there is an algebraic subvariety $Z'\subset Z$ so that
the map $(z,g)\mapsto gz$ from  $Z'\times G^+$ to $Z$
is an isomorphism,
where $G^+\subset G$ is the subgroup of orientation preserving
elements of $G$.
Then  $\ul$ has $n-1$ fixed vertices,
and these vertices are fixed at points in the fixed subspace of $G$.
\item  Otherwise,  $\ul$ has only $n+1$ fixed vertices,
which we may take to be fixed at 0 and $e_i$, $i=1,\ldots ,n$.
\end{enumerate}
Moreover, there is  an analytic function
$\alpha \colon Z\times F\to ({\mathbb  R}^n)^m$ so that 
$$\conf(\ul)=\{\, (x,\alpha (x,f)) 
\mid x\in Z, \mand\ f\in F\,\}$$
and so that the map $(x,f)\mapsto (x,\alpha (x,f))$ is an analytic isomorphism.

What is more, there is a converse to part 2 above if $n=2$.
If $\ul$ is a connected planar classical linkage with $1$ fixed vertex and at least one
other vertex, then there is an algebraic subset $Z'\subset \conf(\ul)$
and a subgroup $G^+\subset \Eu(2)$ conjugate to $SO(2)$
so that the diagonal action $Z'\times G^+\to \conf(\ul)$ is an isomorphism.
\end{thm}

\begin{rem}
For $n>2$ and $G=O(2)$, the characterization of configuration spaces
must be more complicated than that suggested by part 2 of Theorem \ref{thm:5}.
Consider for example the linkage in ${\mathbb R}^3$ with three vertices
$A$, $B$, and $C$, and edge $AC$ of length 1, and with $A$ fixed at $0$
and $B$ fixed at $e_1$.
Then $\conf(\ul)=0\times e_1\times S^2$ which does not satisfy the
hypotheses of part 2 of Theorem \ref{thm:5}.
\end{rem}

In the planar case $n=2$ we may complete our description of
configuration spaces by extending Theorem \ref{thm:5} to the noncompact case.
We let $\Eu(n)^+$ denote the subgroup of $\Eu(n)$ consisting
of orientation preserving Euclidean motions.

\begin{thm}\label{thm:5.2}
Suppose $Z\subset ({\mathbb R}^2)^k$ is an algebraic set invariant under the action of $\Eu(2)$.
Suppose there is a compact algebraic subset $Z'\subset Z$ so that the map
$\beta \colon Z'\times \Eu(2)^+\to Z$ is an isomorphism where
$\beta (z,g)=gz$.
Then there is a connected classical linkage $\ul$ with no fixed vertices
and a finite set $F$ so that
$\conf(\ul)$ is analytically isomorphic to $Z\times F$.
Moreover, there is  an analytic function
$\alpha \colon Z\times F\to ({\mathbb  R}^2)^m$ so that 
$$\conf(\ul)=\{\, (x,\alpha (x,f)) 
\mid x\in Z, \mand\ f\in F\,\}$$
and so that the map $(x,f)\mapsto (x,\alpha (x,f))$ is an analytic isomorphism.

Conversely, if $\ul$ is a connected classical planar linkage with no fixed vertices
and at least two vertices
then there is a compact algebraic subset $Z'\subset \conf(\ul)$ so that
the map $(z,g)\mapsto gz$ gives an isomorphism from
$Z'\times \Eu(2)^+$ to $\conf(\ul)$.
\end{thm}

\section{Functional Linkages}

An essential ingredient in the proofs of the above theorems
is the notion of a functional linkage.
A (quasi)functional linkage is a linkage which ``computes'' some function.

A linkage $\ul$ is \fdf{\qf\ } for a map $f\colon ({\mathbb  R}^n)^k\to ({\mathbb  R}^n)^m$ 
 if there are 
 vertices $w_1,\ldots ,w_k$ and $v_1,\ldots ,v_m$ of $\ul$ so that
 if $p\colon \conf(\ul)\to ({\mathbb  R}^n)^m$ is 
 $p(\varphi )=(\varphi (v_1),\ldots ,\varphi (v_m))$ and
 $q\colon \conf(\ul)\to ({\mathbb  R}^n)^k$ is 
$q(\varphi )=(\varphi (w_1),\ldots ,\varphi (w_k))$ 
then $p=f\compose q$.
The set $q(\conf(\ul))$ is called the \fdf{domain} of the \qf\ linkage.
We call $q$  the input map and call $p$ the output map.

If in addition, there is a $U\subset q(\conf(\ul))$ 
so that the restriction $q|\colon q^{-1}(U)\to U$ is an analytically trivial
 covering map, we say that 
$\ul$ is \fdf{functional} for $f$ 
 with \fdf{\rd\ } $U$.

Moreover, if $q\colon \conf(\ul)\to q(\conf(\ul))$ is an analytic isomorphism
we say that $\ul$ is \fdf{strongly functional} for $f$.
In this case, for expository convenience, if
 $U\subset q(\conf(\ul))$ we say that $\ul$ is  strongly functional for $f$ 
  with \rd\ $U$.

 We call $w_1,\ldots ,w_k$ the \fdf{input vertices} and call
 $v_1,\ldots ,v_m$ the \fdf{output vertices}.
 Repetitions of vertices are allowed,
 although they are not necessary for the results in this paper.

So if $\ul$ is functional for $f$, then 
over $U$ the configuration space is a bunch of copies of the graph
 of $f$.
If the configuration space is just one copy of the graph
 of $f$ it is strongly functional.

The following is a key to the proofs of the above theorems.
Its proof will occupy a substantial part of this paper.

\begin{thm}\label{thm:6}
Suppose $f\colon ({\mathbb  R}^n)^k\to ({\mathbb  R}^n)^m$ is a polynomial map and 
$K\subset ({\mathbb  R}^n)^k$ is compact, $n\ge 2$.
Then there is a functional classical linkage $\ul$ for $f$ with \rd\
$K$.
There is also a strong functional \rl\ $\ul'$ for $f$ with \rd\
$K$.
We may specify that all input and output vertices of the
functional linkages $\ul$ and $\ul'$ be distinct.
\end{thm}

\section{Functoriality of $\conf(\ul)$ and $\sconf(\ul,W)$}

Let $\ul'\subset \ul$ be a sublinkage.
This means that $L'\subset L$, $\ell '=\ell |_{\edge(L')}$, 
$V'\subset V$, $\mu '=\mu |_{V'}$,
and $F'=F\cap \edge(L')$.
Then we have a natural map $\rh{}{\ul'}\colon \conf(\ul)\to \conf(\ul')$
obtained by restriction, i.e., $\rh{}{\ul'}(\varphi )=\varphi |_{\vrtcs(L')}$.
If $L'$ is a single vertex $v$ of $\ul$ and has no edges and $V'$ is empty, then
 we denote $\rh{}{\ul'}=\rh{}v$.
 Thus $\rh{}v(\varphi )=\varphi (v)\in {\mathbb  R}^n=\conf(\ul')$.

If $\ul'\subset \ul$ and $\ul''\subset \ul$ are two sublinkages then we
may define their union $\ul'\cup \ul''$
as the sublinkage $(L''',\ell ''', V''',\mu ''',F''')$ of $\ul$
with $L'''=L'\cup L''$ and $V'''=V'\cup V''$.
Similarly, we may define the intersection $\ul'\cap \ul''$.

\begin{lem}\label{lem:6.1}
If $\ul'\subset \ul$ and $\ul''\subset \ul$ are two sublinkages then we have
a natural identification of $\conf(\ul'\cup \ul'')$ with the fiber product
of the restriction maps $\rh{'}{\ul'\cap \ul''}\colon \conf(\ul')\to \conf(\ul'\cap \ul'')$
and $\rh{''}{\ul'\cap \ul''}\colon \conf(\ul'')\to \conf(\ul'\cap \ul'')$.
\end{lem}
$$
\begin{CD}
\conf(\ul'\cup \ul'') @>{\quad\quad\quad\ }>> \conf(\ul'') \\
@VVV   @VV{\rh{''}{\ul'\cap \ul''}}V\\
\conf(\ul')     @>{\rh{'}{\ul'\cap \ul''}}>> \conf(\ul'\cap \ul'')
\end{CD}
$$

\begin{proof}
This is because a  realization of $\ul'\cup \ul''$ is just
a  realization of $\ul'$ and a  realization of $\ul''$
which happen to agree on $\ul'\cap \ul''$.
Thus 
\begin{equation}\label{eqn:3}
\conf(\ul'\cup \ul'')=\{(\varphi ',\varphi '')\in \conf(\ul')\times \conf(\ul'') \mid 
\rh{'}{\ul'\cap \ul''}(\varphi  ')=\rh{'}{\ul'\cap \ul''}( \varphi '')\}
\end{equation}
is the fiber product.
Strictly speaking, rather than equality in equation (\ref{eqn:3}) above,
the map $(\rho _{\ul'\cup \ul'',\ul'}\ ,\  \rho _{\ul'\cup \ul'',\ul''})$
gives an isomorphism 
between the the two sides of equation (\ref{eqn:3}).
But we will suppress such distinctions.
\end{proof}

As a consequence of Lemma \ref{lem:6.1}, the (semi)configuration space of
the disjoint union of linkages is the product of their
(semi)configuration spaces.

\begin{lem}\label{lem:6.2}
If $\ul$ is the disjoint union of sublinkages $\ul_i$, $i=1,\ldots ,m$,
then 
\begin{eqnarray*}
\conf(\ul)&=&\prod _{i=1}^m \conf(\ul_i)\\
\sconf(\ul,W)&=&\prod _{i=1}^m \sconf(\ul_i,W\cap \vrtcs(\ul_i))
\end{eqnarray*}
\end{lem}

\begin{lem}\label{lem:6.3}
Let $\ul$ be a \rl\ and let $G\subset \Eu(n)$ be a
subgroup which fixes the images of all fixed vertices of $\ul$.
Then $\conf(\ul)$ and $\sconf(\ul,W)$ are invariant
under the diagonal action of $G$.
In particular, if $\ul$ has no fixed vertices,
then $\conf(\ul)$ and $\sconf(\ul,W)$ are invariant under the action of $\Eu(n)$.
\end{lem}

\begin{proof}
If $\beta \in G$, then $\beta (\ul)=\ul$.
So $\beta (\conf(\ul))=\conf(\beta (\ul))=\conf(\ul)$ and
$\beta (\sconf(\ul,W))=\sconf(\beta (\ul),W)=\sconf(\ul,W)$.
\end{proof}

\begin{lem}\label{lem:3.4}
Suppose $\ul$ is a linkage with no fixed vertices.
Form $\ul'$ from $\ul$ by fixing one of the vertices of $\ul$
to some point.
Then there is an isomorphism $\eta \colon \conf(\ul')\times {\mathbb R}^n\to \conf(\ul)$
where $\eta (\varphi ,z)(v)=\varphi (v)+z$ for any vertex $v$ of $\ul$.
\end{lem}

\begin{proof}
Translation preserves all lengths, so $\eta (\varphi ,z)\in \conf(\ul)$.
But if $\ul'$ is obtained by fixing a vertex $v_0$ to $z_0$,
then the inverse of $\eta $ is given by $\eta ^{-1}(\varphi )=(\psi ,z)$
where $\psi (v)=\varphi (v)-\varphi (v_0)+z_0$ and $z=\varphi (v_0)-z_0$.
\end{proof}

\begin{lem}\label{lem:6.4}
Let $\ul$ be a \rl\ and let
$v_1,\ldots ,v_m$ be vertices of $\ul$ which are not fixed.
Let $\ul'$ be obtained from $\ul$ by fixing the vertices $v_1,\ldots ,v_m$
to be at the points $z_1,\ldots ,z_m$.
Let $p\colon \conf(\ul)\to ({\mathbb  R}^n)^m$ be the map $(\rh{}{v_1},\ldots ,\rh{}{v_m})$.
Then $$\conf(\ul')=p^{-1}(z_1,\ldots ,z_m)$$
If $W\subset \vrtcs(\ul)$ and $v_i\in W$ for all $i$,  and
$q\colon \sconf(\ul,W)\to ({\mathbb  R}^n)^m$ is the map 
$q(\varphi )=(\varphi (v_1),\ldots ,\varphi (v_m))$.
Then $$\sconf(\ul',W)=q^{-1}(z_1,\ldots ,z_m)$$
\end{lem}

\begin{proof}
Note that $\conf(\ul')\subset \conf(\ul)$ and it must be exactly those
$\varphi $ with $\rh{}{v_i}(\varphi )=\varphi (v_i)=z_i$.
The lemma follows.
\end{proof}

\begin{lem}\label{lem:6.5}
Let $\ul$ be a \rl\ and let
$v_1,\ldots ,v_m$ be vertices of $\ul$.
Let $\ul'$ be obtained from $\ul$ by 
adding new vertices $u_1,\ldots ,u_m$ and new flexible edges
$\lineseg{u_i}{v_i}$ of length $b_i$.
Fix the vertices $u_i$ to points $z_i\in {\mathbb  R}^n$.
Let $p\colon \conf(\ul)\to ({\mathbb  R}^n)^m$ be the map $(\rh{}{v_1},\ldots ,\rh{}{v_m})$.
Then $\conf(\ul')$
is isomorphic to 
$$p^{-1}(\{w\in ({\mathbb  R}^n)^m \mid \ b_i\ge |w_i-z_i|,\  i=1,\ldots ,m\})$$
\end{lem}

\begin{proof}
The inclusion $\ul\subset \ul'$ gives the map 
$\rh{'}{\ul} \colon \conf(\ul')\to \conf(\ul)$.
Let $$Y=p^{-1}(\{w\in ({\mathbb  R}^n)^m \mid \ b_i\ge |w_i-z_i|,\  i=1,\ldots ,m\})$$
We have a map $\beta \colon Y\to \conf(\ul')$ defined by
$\beta (\varphi )(v)=\varphi (v)$ for $v$ a vertex of $\ul$ and
$\beta (\varphi )(u_i)=z_i$.
Note that $\rh{'}{\ul} (\conf(\ul'))\subset Y$ and $\beta $ is the inverse of
$\rh{'}{\ul} \colon \conf(\ul')\to Y$.
\end{proof}

The following two lemmas are immediate from the definitions.

\begin{lem}\label{lem:6.6}
If $\ul'\subset \ul$ is a  sublinkage then
the map $\rh{}{\ul'}\colon \conf(\ul)\to \conf(\ul')$ is an (analytic)
isomorphism if and only if it is onto and the position
$\varphi (v)$ of each vertex $v$ of $\ul$ is a rational (resp.~analytic) function
of the positions $\varphi (w_i)$ of the vertices $w_i$ in $\ul'$.
More generally, if $Z\subset \conf(\ul)$ then the restriction
$\rh{}{\ul'}|_Z\colon Z\to \rh{}{\ul'}(Z)$ is an (analytic) isomorphism if and only if 
for $\varphi \in Z$, the position
$\varphi (v)$ of each vertex $v$ of $\ul$ is a rational (resp.~analytic) function
of the positions $\varphi (w_i)$ of the vertices $w_i$ in $\ul'$.
\end{lem}

\begin{lem}\label{lem:6.7}
Let $\ul$ be a linkage and suppose $v$ and $w$ are two vertices of $\ul$.
Suppose that whenever there are edges $\lineseg vu$ and $\lineseg wu$ to the same vertex $u$,
 that $\ell (\lineseg vu)=\ell (\lineseg vw)$.
Suppose also that there is no edge $\lineseg vw$.
Then we may form a linkage $\ul'$ from $\ul$ 
by identifying the vertices $v$ and $w$,
and identifying any edges $\lineseg vu$ and $\lineseg wu$.
Moreover there is a natural identification of $\conf(\ul')$
with $\{\varphi \in \conf(\ul)\mid \varphi (v)=\varphi (w) \}$.
\end{lem}

Note in Lemma \ref{lem:6.7} that if 
$\{\varphi \in \conf(\ul)\mid \varphi (v)=\varphi (w) \}$
is nonempty, then the hypotheses of Lemma \ref{lem:6.7}
must be satisfied.

\begin{lem}\label{lem:6.13}
If $\sconf(\ul,W)$ is compact and nonempty, then there is a path
in $L$ from each point of $W$ to a  vertex which is fixed.
\end{lem}

\begin{proof}
By Lemma \ref{lem:6.2}, it suffices to assume that $\ul$ is connected.
If $\ul$ has no fixed vertices, we may translate any 
realization to get another realization.  So $\sconf(\ul,W)$
would be noncompact, c.f.~Lemma \ref{lem:3.4}.
So $\ul$ must have a fixed vertex, and by connectedness
there is a path from each vertex of $W$ to this fixed vertex.
\end{proof}

\section{Constructing Polynomial Functional Linkages}

Simple modifications of two functional linkages
allow us to construct a functional linkage for
their composition or cartesian product.

\begin{lem}\label{lem:6.8}
Let $\ul$ and $\ul'$ be (strong) functional linkages
for functions $f\colon ({\mathbb  R}^n)^k\to ({\mathbb  R}^n)^m$ and 
$g\colon ({\mathbb  R}^n)^m\to ({\mathbb  R}^n)^\ell $
with \rd s $U$ and $U'$.
Suppose that $U\cap f^{-1}(U')$ is nonempty.
We may form a linkage $\ul''$ by taking the disjoint union
of $\ul$ and $\ul'$ and then identifying each
output vertex of $\ul$ with the corresponding
input vertex of $\ul'$.
Then $\ul''$ is a (strong) functional linkage
for $g\compose f$ with \rd\ $U\cap f^{-1}(U')$.
\end{lem}

\begin{proof}
Let $v_1,\ldots ,v_m$ be the input vertices of $\ul'$ and let
$w_1,\ldots ,w_m$ be the output vertices of $\ul$.
Suppose first that any duplications in the $v_i$ correspond to
duplications in the $w_i$ and vice versa.
So $v_i=v_j$ if and only if $w_i=w_j$.
Then $\ul''$ is the union of 
$\ul$ and $\ul'$, and their intersection is 
the linkage with no edges and with vertices $V=\{v_1,v_2,\ldots ,v_m\}$.
Let $\rho _1$ and $\rho _2$ be the input and output
maps of $\ul$ and let $\rho _3$ and $\rho _4$ be the
input and output maps of $\ul'$.
By Lemma \ref{lem:6.1}, we know that $\conf(\ul'')$ is the fiber product of
$\rh{}{\ul\cap \ul'}$ and $\rh{'}{\ul\cap \ul'}$.
If $q\colon \conf(\ul\cap \ul')\to ({\mathbb  R}^n)^m$ is the map 
$q=(\rho _{\ul\cap \ul',v_1}\ ,\ldots ,\ \rho _{\ul\cap \ul',v_m})$,
then $\rho _2=q\compose \rh{}{\ul\cap \ul'}$ and $\rho _3=q\compose \rh{'}{\ul\cap \ul'}$.
If $v_1,\ldots ,v_m$ are all distinct, $q$ will be the identity, but if there
are duplications, $q$ will be some sort of diagonal map.
Since $q$ is injective, $\conf(\ul'')$ is also 
the fiber product of $\rho _2$ and $\rho _3$,
\begin{equation}\label{eqn:6}
\conf(\ul'')=\{(\varphi ,\varphi ')\in \conf(\ul)\times \conf(\ul') 
\mid \rho _2(\varphi )=\rho _3(\varphi ')\}
\end{equation}
so that $\rh{''}{\ul}$ and $\rh{''}{\ul'}$ are induced by projection.
Note that 
$$g\compose f\compose \rho _1\compose \rh{''}\ul
=g\compose \rho _2\compose \rh{''}\ul
=g\compose \rho _3\compose \rh{''}{\ul'}=\rho _4\compose \rh{''}{\ul'}$$
so $\ul''$ is \qf\ for $g\compose f$.
The input map is $\rho _1\compose \rh{''}\ul$ and the output map is
$\rho _4\compose \rh{''}{\ul'}$.

Now let us see that we can take the \rd\ to be $U\cap f^{-1}(U')$.
The restriction of $\rho _1$ to $\rho _1^{-1}(U\cap f^{-1}(U'))$ is an analytically
trivial cover since the restriction to $\rho _1^{-1}(U)$ is, so we only need
show that $\rh{''}\ul$ restricts to an analytically trivial
cover of $\rho _1^{-1}(U\cap f^{-1}(U'))=\rho _1^{-1}(U)\cap \rho _2^{-1}(U')$.
We know that there is a finite set $F$ and an analytic isomorphism
$\sigma  \colon U'\times F\to \rho _3^{-1}(U')$ 
so that $\rho _3\compose \sigma $ 
is projection to $U'$.
Now by (\ref{eqn:6}), we have
\begin{eqnarray*}
\rh{''}\ul^{-1}(\rho _2^{-1}(U'))&=&\{(\varphi ,\varphi ')\mid \rho _2(\varphi )
=\rho _3(\varphi ')\in U' \}\\
&=&\{(\varphi ,\sigma  (\rho _2(\varphi ),c))\mid \rho _2(\varphi )\in U' 
\mand c\in F\}
\end{eqnarray*}
So we have an analytic trivialization 
$\sigma  '\colon \rho _2^{-1}(U')\times F\to \rh{''}\ul^{-1}(\rho _2^{-1}(U'))$
given by $\sigma  '(\varphi ,c)=(\varphi ,\sigma  (\rho _2(\varphi ),c))$.
So $\ul''$ is functional for $g\compose f$ with \rd\ $U\cap f^{-1}(U')$.

To prove  the strong case, note that all the covers are one-fold and hence
are analytic isomorphisms
and we may take $U=\rho _1(\conf(\ul))$ and $U'=\rho _3(\conf(\ul'))$.

If there are duplications in the input and output vertices things can get
more complicated,
since we may end up having to identify vertices in $\ul$ or $\ul'$
which were not previously identified.
Let $\Delta _{ij}=\{ (z_1,\ldots ,z_m)\in ({\mathbb  R}^n)^m \mid z_i=z_j\}$.

Suppose $v_i=v_j$ but $w_i\neq w_j$.
 Then we must have $U'\subset \Delta _{ij}$.
Also, in $\ul''$ we end up identifying $w_i$ with $w_j$.
Let us first see whether we can do so according to Lemma \ref{lem:6.7}.
Suppose $w$ is another vertex so that $\lineseg {w}{w_i}$ and $\lineseg{w}{w_j}$
are both edges of $\ul$.  
Since $U\cap f^{-1}(U')$ is nonempty,
there is a $\varphi \in \conf(\ul)$ so that $\rho _1(\varphi )\in U\cap f^{-1}(U')$.
Hence $\rho _2(\varphi )=f(\rho _1(\varphi ))\in U'\subset \Delta _{ij}$,
 and so $\varphi (w_i)=\varphi (w_j)$.
So 
$$\ell (\lineseg{w}{w_i})=|\varphi (w)-\varphi (w_i)|
=|\varphi (w)-\varphi (w_j)|=\ell (\lineseg{w}{w_j})$$
So in $\ul''$ we may identify the edges $\lineseg w{w_i}$ and $\lineseg{w}{w_j}$
since they have the same length.
There could not be an edge $\lineseg{w_i}{w_j}$ since 
$0\neq \ell (\lineseg{w_i}{w_j})=|\varphi (w_i)-\varphi (w_j)|=0$.
So by Lemma \ref{lem:6.7} we are allowed to take the quotient linkage
$\ul_1$ of $\ul$, identifying $w_i$ and $w_j$.
By Lemma \ref{lem:6.7} we also see that $\ul_1$ is still functional for 
$f$ but the domain has 
shrunk from $\rho _1(\conf(\ul))$ to $\rho _1(\conf(\ul))\cap f^{-1}(\Delta _{ij})$.
So we may take the restricted domain of $\ul_1$ to be $U\cap f^{-1}(\Delta _{ij})$.
Do this identification for each pair $i,j$ with $v_i=v_j$ and $w_i\neq w_j$
and we
eventually get a functional linkage $\ul_2$ for $f$ with restricted domain
$U_2=U\cap f^{-1}(\Delta )$ for some $\Delta \supset U'$.

Now suppose $w_i=w_j$, but $v_i\neq v_j$.
Then we must have $f(U)\subset \Delta _{ij}$.
Also, in $\ul''$ we end up identifying $v_i$ with $v_j$.
Let us see whether we can do so.
Suppose $v$ is another vertex so that $\lineseg{v}{v_i}$ and $\lineseg{v}{v_j}$
are both edges of $\ul'$.  
Since $U\cap f^{-1}(U')$ is nonempty, we know $\Delta _{ij}\cap U'$
is nonempty, so
there is a $\varphi \in \conf(\ul')$ so that $\rho _3(\varphi )\in \Delta _{ij}\cap U'$,
and hence $\varphi (v_i)=\varphi (v_j)$.
So as above, Lemma \ref{lem:6.7} will allow us to take the quotient linkage
$\ul'_1$ identifying $v_i$ and $v_j$.
By Lemma \ref{lem:6.7} we also see that $\ul'_1$ is still functional for 
$g$ but with restricted domain  $U'\cap \Delta _{ij}$.
Do this identification for each pair $i,j$ with $w_i=w_j$ and $v_i\neq v_j$
and we
eventually get a functional linkage $\ul'_2$ for $g$ with restricted domain
$U'_2=U'\cap \Delta '$ for some $\Delta '\supset f(U)$.

After doing all these identifications, we have $\ul''$ is the union of $\ul_2$ and $\ul'_2$,
and we may finish the proof as above.
The only thing to check is that $U_2\cap f^{-1}(U'_2)=U\cap f^{-1}(U')$.
But $U_2\cap f^{-1}(U'_2)=U\cap f^{-1}(\Delta )\cap f^{-1}(U')\cap f^{-1}(\Delta ')=U\cap f^{-1}(U')$
since $U'\subset \Delta $ and $U\subset f^{-1}(\Delta ')$.
\end{proof}

\begin{lem}\label{lem:6.9}
For $i=0,1$, let $\ul_i$ be (strong) functional linkages
for  functions $f_i\colon ({\mathbb  R}^n)^{k_i}\to ({\mathbb  R}^n)^{m_i}$
with \rd\ $U_i$.
Form a linkage $\ul$ by taking  the disjoint union
of $\ul_0$  and $\ul_1$.
If $k_0=k_1$, form a linkage $\ul'$ by taking $\ul$ and identifying corresponding input
vertices of $\ul_0$ and $\ul_1$.
\begin{enumerate}
 \item
 Then $\ul$ is a (strong) functional linkage
for $$f_0\times f_1\colon ({\mathbb  R}^n)^{k_0}\times ({\mathbb  R}^n)^{k_1}\to 
({\mathbb  R}^n)^{m_0}\times ({\mathbb  R}^n)^{m_1}$$ with \rd\ $U_0\times U_1$.
(In particular, if $m_1=0$ then $f_0\times f_1$ is the composition of $f_0$ with
projection 
$({\mathbb  R}^n)^{k_0}\times ({\mathbb  R}^n)^{k_1}\to ({\mathbb  R}^n)^{k_0}$.)
\item
If $k_0=k_1$, $\ul'$ is a (strong) functional linkage
for $$(f_0,f_1)\colon ({\mathbb  R}^n)^{k_0}\to 
({\mathbb  R}^n)^{m_0}\times ({\mathbb  R}^n)^{m_1}$$ with \rd\ $U_0\cap U_1$.
\end{enumerate}
\end{lem}

\begin{proof}
The statement for $\ul$ is trivial to prove.
If $k_0=k_1=k$, let $\ul''$ be the linkage with no edges 
and $k$ vertices $v_1,\ldots ,v_k$.
If we let these vertices be the input vertices and let the output vertices
be doubled, $v_1,\ldots ,v_k,v_1,\ldots ,v_k$, then we get a strong functional
linkage for the diagonal map $\Delta (z)= (z,z)$
with domain all of $({\mathbb R}^n)^k$.
By Lemma \ref{lem:6.8}, $\ul'$ is the composition of the linkages $\ul$ and $\ul''$
and hence
is (strongly) functional for 
$(f_0\times f_1)\compose \Delta =(f_0,f_1)$.
\end{proof}

\subsection{Elementary Polynomial Functional Linkages}

We are now ready to make the first progress in proving Theorem \ref{thm:6}.
We first reduce it to finding (strong) functional linkages for
addition, multiplication, and some linear maps.

\begin{reduc}\label{red:1}
To prove Theorem \ref{thm:6}, it suffices to prove the existence of
(strong) functional linkages for the following functions,
each with arbitrarily large compact \rd,
and distinct input and output vertices.
\begin{enumerate}
\item  $q\colon ({\mathbb R}^n)^2\to {\mathbb R}^n$ given by $q(x,y)=x+y$.
\item  $r\colon T\times T\to T$ given by $r(sz_0,tz_0)=stz_0$, where $z_0\neq 0$ and
   $T$ is the line through $0$ and $z_0$.
\item  $u\colon {\mathbb R}^n\to {\mathbb R}^n$ 
any rank 1 linear transformation.
\end{enumerate}
\end{reduc}

\begin{proof} 
Suppose $f\colon ({\mathbb R}^n)^k\to ({\mathbb R}^n)^m$ is a polynomial map
for which we wish to find a (strong) functional linkage.
By Lemma \ref{lem:6.9},
it suffices to consider the case $m=1$, i.e.,
of polynomials $f\colon ({\mathbb R}^n)^k\to {\mathbb R}^n$.
Note that if $\ul_i$ are functional linkages for
$f_i\colon ({\mathbb  R}^n)^k\to {\mathbb  R}^n$, by Lemma \ref{lem:6.9}
 we may form a functional linkage
for $(f_0,f_1)\colon ({\mathbb  R}^n)^k\to {\mathbb  R}^n\times {\mathbb  R}^n$, 
then using Lemma \ref{lem:6.8} and composing with a functional linkage for
$q(x,y)=x+y$ we get a functional linkage for 
$f_0+f_1\colon ({\mathbb  R}^n)^k\to {\mathbb  R}^n$.
So it suffices to find
functional linkages for 
$f\colon ({\mathbb  R}^n)^k\to {\mathbb  R}^n$ of the form $f(x)=p(x)e_j$
for $p$ a monomial.
If $\ul_i$ are functional linkages for
$f_i\colon ({\mathbb  R}^n)^k\to {\mathbb  R}^n$, of the form $f_i(x)=p_i(x)e_j$,
then by Lemma \ref{lem:6.9}, Lemma \ref{lem:6.8} and
the map $r$ above, we get a functional linkage for the map
$x\mapsto p_0(x)p_1(x)e_j$.
So it suffices to find functional linkages for degree $0$ or $1$ monomials.
Degree 1 monomials are linear and so are obtained from the map $u$ above.
In particular, we take a functional linkage for $u$ and add
$k-1$ disjoint vertices which we designate as input vertices.
Degree 0 monomials are constants, use a trivial linkage with one fixed
output vertex and $k$ input vertices.
\end{proof}

While the functions we reduced to in Reduction \ref{red:1}
seem natural, they are not all suited to easy description as linkages.
So we make a further reduction to some elementary functions
for which we can more readily provide linkages.
It is interesting to note that by  item 4 below, we could actually construct functional
linkages for any rational function.
Later in this paper we will investigate further exactly
which functions admit functional linkages.

\begin{reduc}\label{red:2}
To prove Theorem \ref{thm:6}, it suffices to prove the existence of
(strong) functional linkages for the following functions,
all with distinct input and output vertices.
\begin{enumerate}
\item  Translation:  $z\mapsto z+z_0$, with \rd\ any compact $K\subset {\mathbb  R}^n$.

\item  Scalar multiplication: $z\mapsto \lambda z$, with \rd\ a disc
  $\{ |z-z_0|\le r \}$ for some $z_0$ and for any $r$ as large as we wish.
  
 \item  Average:  $(z,w)\mapsto (z+w)/2$, with \rd\
 $\{(z,w)\mid |z-z_0|\le r, |w+z_0|\le r\}$  for some $z_0$ and for any 
 $r$ as large as we wish.
 
 \item  Inversion in a line:  $sz_0\mapsto (1/s) z_0$, for any specified  $z_0\in {\mathbb  R}^n-0$,
with \rd\
 any compact $K\subset \{tz_0 \mid t\neq 0\}$.
 

\item  Orthogonal projection to a line: $z\mapsto (z\cdot z_0)z_0$, 
for any unit vector $z_0\in {\mathbb  R}^n$, with \rd\ any compact $K$.
\end{enumerate}
\end{reduc}

\begin{proof}
Note first that we may always further restrict the domain of a
functional linkage, so it suffices to find functional linkages with
arbitrarily large compact \rd s, for example
(products of) balls of radius $r$.
By Reduction \ref{red:1}, we only need to use the above
five types of functional linkages to construct functional linkages
for the three types
of functions listed in Reduction \ref{red:1}.

But before we do this, we will show that for scalar multiplication 2 above,
we may actually take the \rd\ to be an arbitrarily large ball
$\{ |z|\le r\}$.
First
use 2 above to get a functional linkage for $z\mapsto \lambda z$, with \rd\
  $\{ |z-z_0|\le r \}$ for some $z_0$.  Then use 1 to get a functional linkage for
 translation  $z\mapsto z+z_0$ with \rd\ $|z|\le r$.
  Using Lemma \ref{lem:6.8}, compose these two to get a functional
  linkage for $z\mapsto \lambda z+\lambda z_0$ with \rd\ $|z|\le r$.
  Now using 1 and Lemma \ref{lem:6.8}, compose with a translation by $-\lambda z_0$
  to get our desired linkage for $z\mapsto \lambda z$ with \rd\ $|z|\le r$.

To get $(z,w)\mapsto z+w$ with \rd\ 
$|z|\le r$, $|w|\le r$, find a functional
linkage for the average 3 above,
with \rd\
 $\{(z,w)\mid |z-z_0|\le r, |w+z_0|\le r\}$  for some $z_0$.
  Then using 1, find functional linkages
for $z\mapsto z+z_0$ and $z\mapsto z-z_0$, both with \rd\ $|z|\le r$.
By Lemma \ref{lem:6.9}, their disjoint union is functional for
$(z,w)\mapsto (z+z_0,w-z_0)$ with \rd\ $|z|\le r$, $|w|\le r$.
Using Lemma \ref{lem:6.8} and composing with the first linkage, we get
a functional linkage for $(z,w)\mapsto (z+w)/2$ with \rd\ $|z|\le r$, $|w|\le r$.
Now composing with scalar multiplication by $2$ with \rd\ $|z|\le r$,
 we get a functional
linkage for $(z,w)\mapsto z+w$ with \rd\ $|z|\le r$, $|w|\le r$.

Next we will find a (strong) functional linkage for any rank one linear
map $u\colon {\mathbb R}^n\to {\mathbb R}^n$.
But any rank one linear map is a composition of orthogonal
projections to lines followed by scalar multiplication.
So compositions of maps 5 and 2 will give us a (strong)
functional linkage for $u$.
(To see that $u$ is such a composition, first do orthogonal 
projection to the orthogonal compliment of $\ker u$.  
If $u^2\neq 0$, we may then orthogonal
project to the image of $u$, and multiply by an appropriate
scalar to get $u$.  If $u^2=0$, do two more orthogonal projections,
first to a line $T$ which is neither perpendicular to the image of $u$
nor contained in $\ker u$, and then orthogonal project to the
image of $u$. Finally, multiply by an appropriate scalar.)

So the only remaining function is multiplication.
We will do some algebraic manipulation to get multiplication
$(sz_0,tz_0)\mapsto stz_0$.
First, note that $$st=((s+t)^2-(s-t)^2)/4$$
So it suffices to find a functional linkage for $sz_0\mapsto s^2z_0$ with
\rd\ 
$\{ \,  |s|\le r\,\}$.
By 4 above there is a functional linkage for the map
$h(sz'_0)=(1/s)z'_0$ with \rd\ $2/3\le s\le 4/3$,
where $z_0'=3rz_0$.
But we note that 
$$z_0'-h((h((1+s)z_0')+h((1-s)z_0'))/2)=s^2z_0'$$
so by composing known functional  linkages,
we get a functional  linkage for the function $sz_0'\mapsto s^2z_0'$
with \rd\ $|s|\le 1/3$.
But this is the same as the function $sz_0\mapsto (s^2/(3r))z_0$
with \rd\ $|s|\le r$.
So after composing with multiplication by $3r$,
we get the desired functional  linkage.
\end{proof}

So we have now reduced 
the proof of Theorem \ref{thm:6} to finding the functional linkages 1-5 
in Reduction \ref{red:2}.
In doing so, the following Lemma will be useful.
Its proof may be safely left to the reader.
It is, for example, a special case of the theorem that
a proper submersion is a locally trivial fibration.

\begin{lem}\label{lem:6.10}
Let $f\colon M\to {\mathbb  R}^n$ be a smooth map from a compact
$n$ dimensional manifold with boundary.
Let $S\subset M$ be the set of critical points of $f$,
the points where $df$ has rank $<n$.
Let $U$ be any connected component of ${\mathbb  R}^n-f(S\cup \partial M)$.
Then the restriction $f|\colon f^{-1}(U)\to U$ is a covering projection.
\end{lem}

In our usage, $f$ is analytic, and $U$
is often contractible, so $f$ restricts to an analytically trivial covering
of $U$, thus $f^{-1}(U)$ is analytically isomorphic to $U\times $ a finite set.
As another application, we will use the consequence that
$f(M)$ is the union of $f(S\cup \partial M)$ and some connected components
of ${\mathbb  R}^n-f(S\cup \partial M)$.

So in the remainder of this section we will construct the functional
linkages 1-5 listed in Reduction \ref{red:2} above.
But first we look at some useful examples.

\subsection{Simulating interior joints, cables, and telescoping edges}

In our model of linkages, edges are connected only at their ends.
Actual linkages used in real life might have a connection in the
middle of an edge.
This may be simulated as in Figure \ref{fig:1}.
If $\ell (\lineseg AB)=b$ and we wish to place a connection $C$
in the middle of $\lineseg AB$,
let $\ell (\lineseg AC)=a$ and $\ell(\lineseg BC)=b-a$.
Thus when drawing linkages, it is allowable to draw
a joint in the middle of an edge.

\begin{figure}
\Myepsf{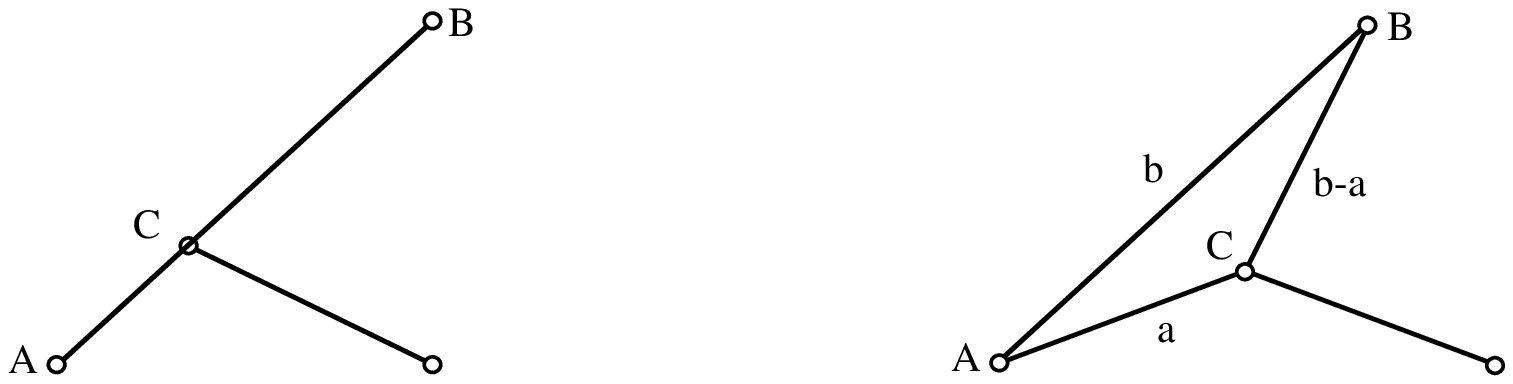 scaled 700}
\capshun 1 { How to put a joint in the middle of an edge}
\end{figure}

If we are in the context of semiconfiguration spaces, we can also
simulate other types of linkages.
For example, suppose we want two vertices $A$ and $B$ connected by a cable, so 
the distance between them is constrained to be $\le b$.
More generally, suppose we wish to connect $A$ and $B$
by a telescoping edge,
so the distance between them is constrained to be
in the interval $[a,b]$.
This can be simulated as in Figure \ref{fig:2}.  Since we are using semiconfiguration
spaces, we can ignore the position of the vertex $D$.
To simulate a cable, we take $c=d=b/2$.
To simulate a telescoping edge with $0<a<b$, we take 
$c=(a+b)/2$, $d=(b-a)/2$.

\begin{figure}
\Myepsf{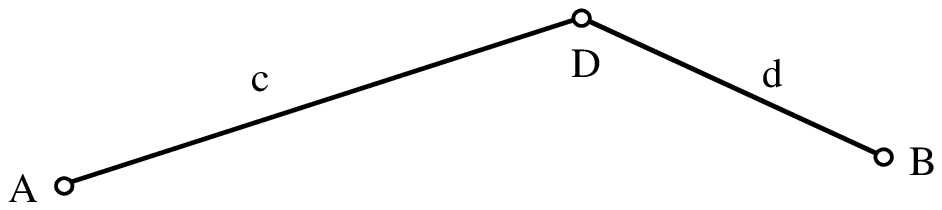 scaled 700}
\capshun 2 { Simulating a cable or telescoping edge}
\end{figure}

\subsection{The rigidified parallelogram}

When constructing linkages, one often wants four vertices to lie
in a plane, and moreover to form a parallelogram.
The linkage of Figure \ref{fig:3} will do this.
There are six vertices $A,B,C,D,E,F$ and nine edges
$\lineseg AB$ and $\lineseg CD$ of length $a$,
$\lineseg AC$, $\lineseg BD$, and $\lineseg EF$ of length $b$
and $\lineseg AE$, $\lineseg BE$, $\lineseg CF$, and $\lineseg DF$
of length $a/2$.

\begin{figure}
\Myepsf{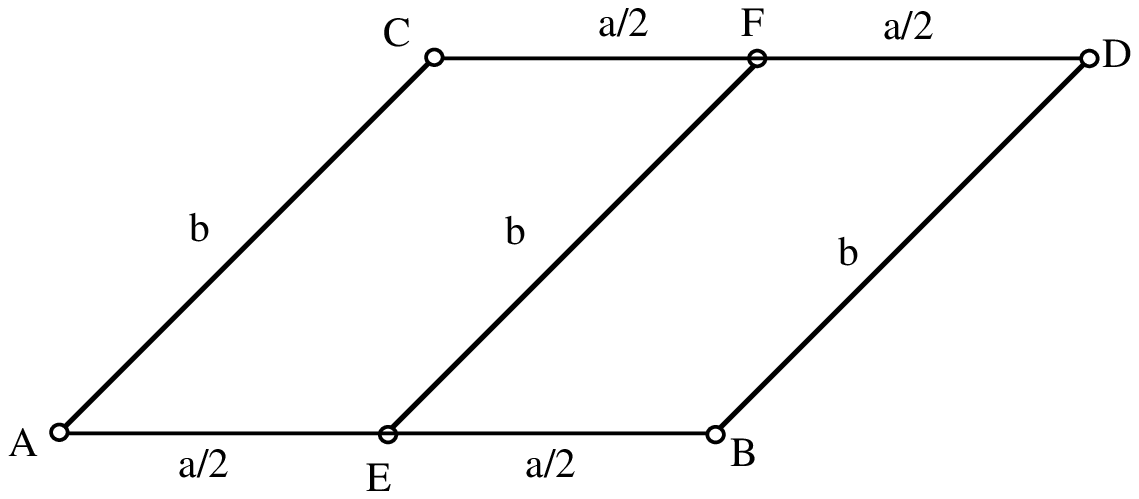 scaled 700}
\capshun 3 { A  rigidified parallelogram}
\end{figure}

Note that $E$ and $F$ are simulated interior joints,
so for any realization $\varphi $,
we must have $\varphi (E)=(\varphi (A)+\varphi (B))/2$ and  
$\varphi (F)=(\varphi (C)+\varphi (D))/2$.

We claim that $\varphi (A),\varphi (B),\varphi (C),\varphi (D)$ form a parallelogram.
To see this, let $x=\varphi (B)-\varphi (A)$, $y=\varphi (C)-\varphi (A)$, and $z=\varphi (D)+\varphi (A)-\varphi (C)-\varphi (B)$.
Then $\varphi (E)=\varphi (A)+x/2$, $\varphi (D)=\varphi (A)+x+y+z$, and $\varphi (F)=\varphi (A)+y+(x+z)/2$.
The side length equations become
$|x|=a$, $|x+z|=a$, $|y|=b$, $|y+z|=b$, and $|y+z/2|=b$.
From the last three equations we see that $z=0$ and hence
that $\varphi (A),\varphi (B),\varphi (C),\varphi (D)$ form a parallelogram.

Henceforth, when drawing such a rigidified parallelogram we will
draw the edge $\lineseg EF$ as a gray line in an attempt to
unclutter the drawings.
We will also usually refrain from naming the vertices $E$ and $F$.

If $b=a$, we will often refer to this
as a rigidified square, (although realizations usually do not have
right angles).
We take this terminology rigidified square or rigidified parallelogram from \cite{KM}.
Note however that the linkage is not completely rigid, but retains some
flexibility.
In fact, the quotient space $\conf(\ul)/Euc(n)$
is an interval, parameterized by the angle at a vertex.
If the rigidifying edge $\lineseg EF$ were not present, the
configuration space would be bigger, including configurations obtained
by bending along $\lineseg CB$ or $\lineseg AD$
which one generally does not want.

\subsection{Making spheres}

\begin{lem}\label{lem:6.11}
Given any round $k$-sphere $S\subset {\mathbb  R}^n$, 
there is a classical linkage $\ul$ and a vertex
$v\in \vrtcs(\ul)$ so that $S=\sconf(\ul,\{v\})$.
Moreover $\rh{}v\colon \conf(\ul)\to {\mathbb  R}^n$ is an isomorphism 
to its image $S$.
\end{lem}

\begin{proof}
After translation and rotation we may assume that 
$$S=\{ x\in {\mathbb  R}^n\mid x_i=0, i=1,\ldots ,n-k-1 \mand |x|=r\,\}$$
Consider a linkage $\ul$ with fixed vertices $v_i$, $i=0,\ldots ,n-k-1$
and one other vertex $v$, and edges $\lineseg v{v_i}$, $i=0,\ldots ,n-k-1$.
We fix $v_0$ at 0 and fix $v_i$ at the point $r e_i$ if $i\ge 1$.
We let the length of $\lineseg v{v_0}$ be $r$, and let the length of
$\lineseg v{v_i}$ be $\sqrt 2 r$ for $i\ge 1$.
Then if $\varphi \in \conf(\ul)$ and $\varphi (v)=x$, we know that
$r=\ell (\lineseg v{v_0})=|x|$ and
$$2r^2=|x-re_i|^2=|x|^2-2rx_i+r^2=2r^2-2rx_i$$
for all $i=1,\ldots ,n-k-1$ so we know that $x_i=0$ for $i\le n-k-1$, so $x\in S$.
Conversely, if $x\in S$ then $x=\varphi (v)$ for some $\varphi \in \conf(\ul)$.

The map $\rh{}v$ is an isomorphism to $S$ since $v$ is the only vertex which is
not fixed.
\end{proof}

\subsection{A simple Linkage, a key to understanding more complicated Linkages}
\label{sec:4.5}

It will be useful to look first at a simple linkage $\ul$,
as shown in the left half of Figure \ref{fig:4}.
Using Lemma \ref{lem:6.11}, we start with a linkage with $n-1$ fixed vertices,
 $A$ and $A_i$, 
$i=3,\ldots , n$, and a movable vertex $B$,
so that $B$ is constrained to move in a circle with center $A$ and radius $a$.
We then add a final vertex $C$ and an edge $\lineseg BC$ of length $b$.
We assume that $b\le a$.

\begin{figure}
\Myepsf{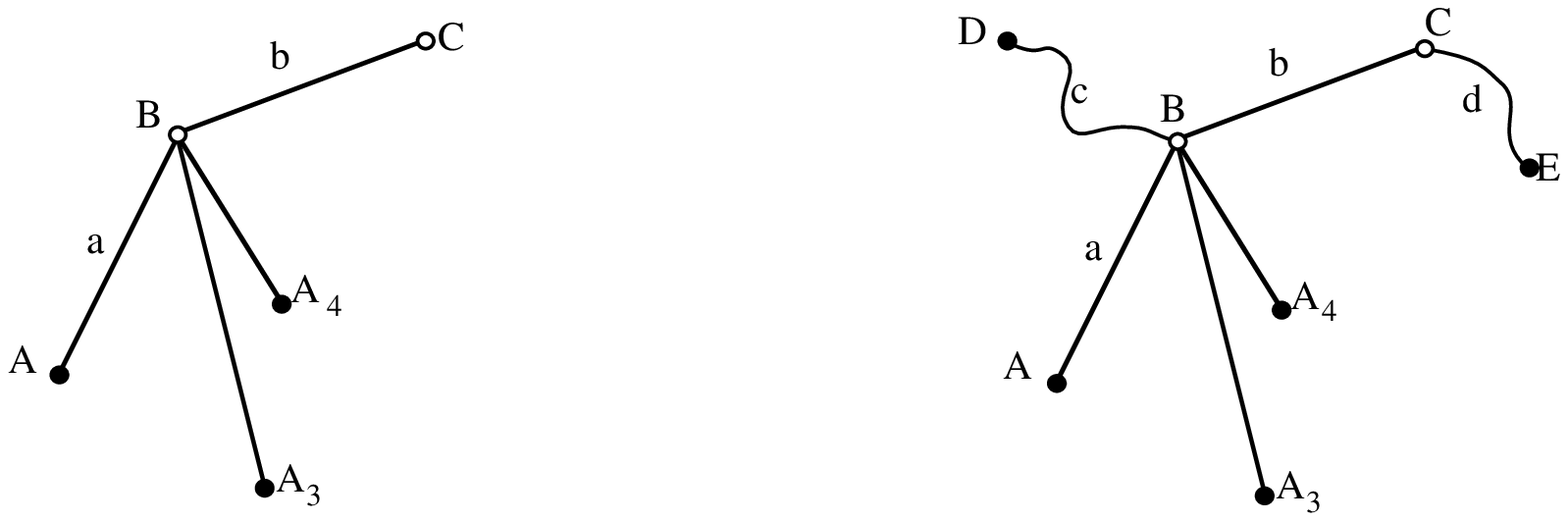 scaled 700}
\capshun 4 { A  Basic Linkage}
\end{figure}

To make some choices, we fix $A$ at a point $z_1$,
fix the $A_i$ vertices at $z_1+a e_i$ for $3\le i\le n$, and put in edges
$\lineseg {A_i}B$ of length  $\sqrt 2 a$ and $\lineseg AB$ of length $a$.
Then $B$ is constrained to lie in a circle 
$$S=\{z\in {\mathbb  R}^n \mid  (z-z_1)\cdot e_i =0 \mfor i\ge 3
\mand  |z-z_1|=a \}$$
and $C$ is only constrained to be in a sphere around $B$ of radius $b$.
Consequently $\conf(\ul)$ is a torus $S^1\times S^{n-1}$,
where $S^1$ is the unit circle about the origin in the $x_1x_2$ plane and
$S^{n-1}$ is the unit sphere about the origin in ${\mathbb R}^n$.
We may identify $(u,v)\in S^1\times S^{n-1}$ with
 $\varphi _{uv}$ where $\varphi _{uv}(A)=z_1$, $\varphi _{uv}(B)=z_1+au$,
 and $\varphi _{uv}(C)=z_1+au+bv$.
Note that $\rh{}C\colon \conf(\ul)\to {\mathbb  R}^n$ is then the map
$\rh{}C(u,v)=z_1+au+bv$ which has critical set
where  $v_1u_2=u_1v_2$.
The image of the critical set is the torus 
of points at distance $b$ from the circle $S$.
So by Lemma \ref{lem:6.10} we see that the image of
$\rh{}C$ is the solid torus  $T$ of points at distance $\le b$ from $S$.
Moreover $\rh{}C$ restricts to a double cover
of the interior of $T$.
In fact this double cover is analytically trivial.
In applications below, we will usually only focus on some disc $\{|z-z_0|\le r\}$
inside the solid torus  where, say, $z_0\in S$
and $0<r<b$.
Then $\rh{}C$ restricts to an analytically trivial
double cover of this disc, (since it is an analytic proper submersion
over the disc, and hence a locally analytically trivial fiber bundle).

When working with \rl s, we will want to modify
this linkage so that $\rh{}C$ is an analytic isomorphism
to some disc $\{|z-z_0|\le d\}$.
We do this by tethering the vertices $B$ and $C$
 to fixed vertices $D$ and $E$ so that their movement
is restricted, see the \rl\ on the right half of Figure \ref{fig:4}.
Consider first the sublinkage $\ul'$ formed by $A$, $A_i$, $B$, $C$,
and $D$, with rigid edges $\lineseg AB$, $\lineseg {A_i}{B}$ and $\lineseg BC$, 
and a flexible
edge $\lineseg BD$ of length $c$, where $D$ is fixed at some point $z_2$
(and $A$ and $A_i$ are fixed as before).
By Lemma \ref{lem:6.5}, we have 
\begin{eqnarray*}
\conf(\ul')&=&\rh{}B^{-1}(\{|z-z_2|\le c\})\\
&=&\{(u,v)\in S^1\times S^{n-1} \mid c \ge |z_1+au-z_2| \}=T\times S^{n-1}
\end{eqnarray*}
for some arc $T$ of $S^1$, as long as we choose $z_2$
and $c$ appropriately.
For convenience, we choose $c=\sqrt 2 a$ and $z_2=z_1+aw_0$
for some $w_0\in S^1$.
Then $T$ will be the semicircle between $\pm w'_0$
which contains $w_0$,
where $w'_0$ is obtained by rotating $w_0$ by $\pi /2$.
By Lemma \ref{lem:6.10},
we know that $\rh{}C$ restricts to an analytically trivial
covering of $\{|z-z_1-aw'_0|<b\}$.
But by checking the inverse image of a point, for example
$z_1+aw'_0$, we see that it is a one-fold cover,
hence an analytic isomorphism.
So now in $\ul$, if we fix $E$ at $z_1+aw'_0$
and pick $d<b$, we see that
$\rh{}C\colon \conf(\ul)\to {\mathbb  R}^n$ is an analytic isomorphism to its
image $|z-z_1-aw'_0|\le d$.

\subsection{A Functional Linkage for Translation }

Now let us find a functional linkage for translation $z\mapsto z+z_0$
with \rd\ $|z|\le r$.
Consider the linkages $\ul$ in Figure \ref{fig:5}, which we will show to be functional
for $z\mapsto z+z_0$ with \rd\ $|z|\le r$.
The right hand \rl\ will be strongly functional.

\begin{figure}
\Myepsf{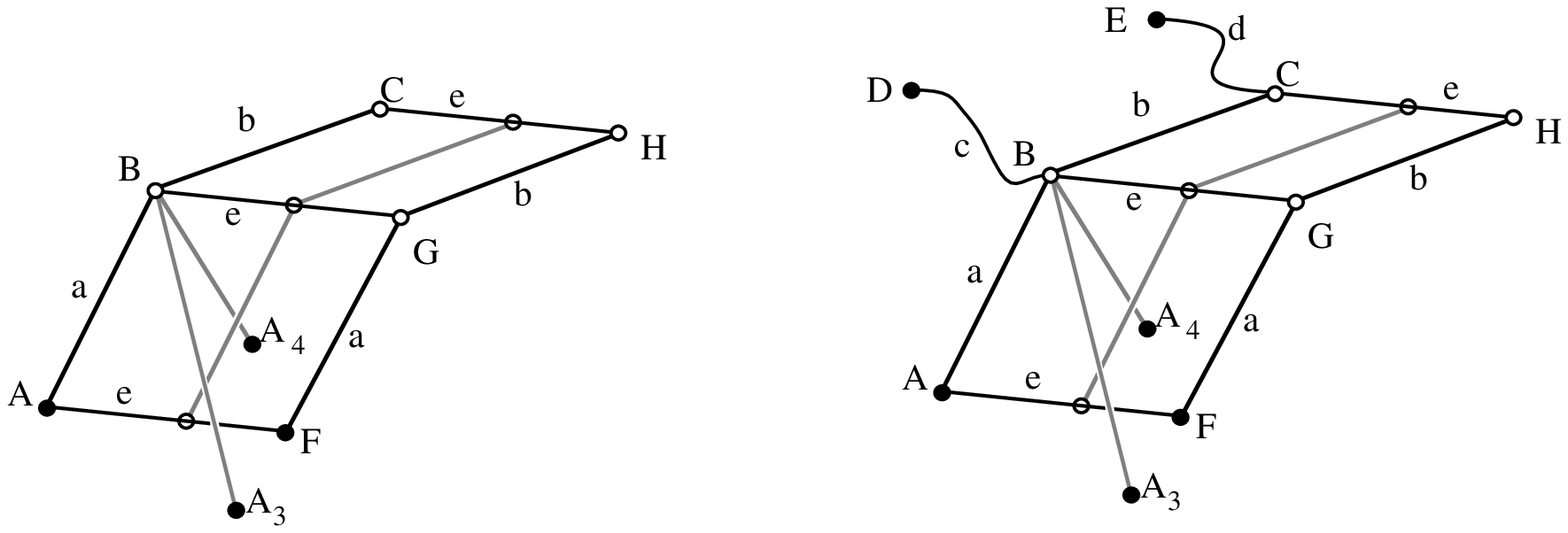 scaled 700}
\capshun 5 {Translation Linkage }
\end{figure}

Choose $a>2r$, $b\le a$, and let $e=|z_0|$.
We start with the sublinkage $\ul'$ which is that of \figsix,
with $z_1$ to be determined later.
We add a vertex $F$ fixed at $z_1+z_0$,
and add vertices $G$ and $H$, and add
edges $\lineseg AF$, $\lineseg BG$ and $\lineseg CH$ of length $e$,
and add edges $\lineseg FG$ and $\lineseg GH$ of lengths $a$ and $b$.
We also rigidify the parallelograms $ABGF$ and $BCHG$.

We let $C$ be the input vertex and $H$ be the output vertex.
Notice $\lineseg AF$, $\lineseg BG$, and $\lineseg CH$ are parallel, and so
for any $\varphi \in \conf(\ul)$ we must have $\varphi (H)=\varphi (C)+z_0$.
So $\ul$ is \qf\ and we must only check that the \rd\ can be $|z|\le r$.

We claim by Lemma \ref{lem:6.6} 
that $\rh{}{\ul'}\colon \conf(\ul)\to \conf(\ul')$ is an isomorphism.
This is because the positions of $F$, $G$, $H$, and the other three
unnamed vertices used to rigidify the quadrilaterals are all
polynomial functions of the positions of $A$, $B$, and $C$.
Now the fact that $\rh{}C$ doubly covers $|z|\le r$ (for the left hand
classical linkage) or singly covers $|z|\le r$ (for the right hand
\rl) follows from the discussion of $\conf(\ul')$
in \figsix,
as long as we make appropriate choices of $z_1$, $a$, $b$,
and $w_0$.
For example, we may choose $b$ so $r<b<a$,
choose $w_0=e_1$, $w_0'=e_2$, and $z_1=-ae_2$.

\subsection{A Functional Linkage for real scalar Multiplication}

Now let us find a functional linkage $\ul$ for scalar multiplication
$z\mapsto \lambda z$.
To do this we marry the pantograph 
of Figure \ref{fig:10} with 
the linkage of \figsix.
To make the pantograph, we take a rigidified parallelogram $EDFB$
with side lengths $b$ and $ca$.
To this we add vertices $A$ and $C$ and sides 
$\lineseg AE$ of length $a$, $\lineseg AD$ of length $a+ca$,
$\lineseg FC$ of length $cb$, and $\lineseg CD$ of length $b+cb$.
Thus $E$ and $F$ are simulated interior joints.

\begin{figure}
\Myepsf{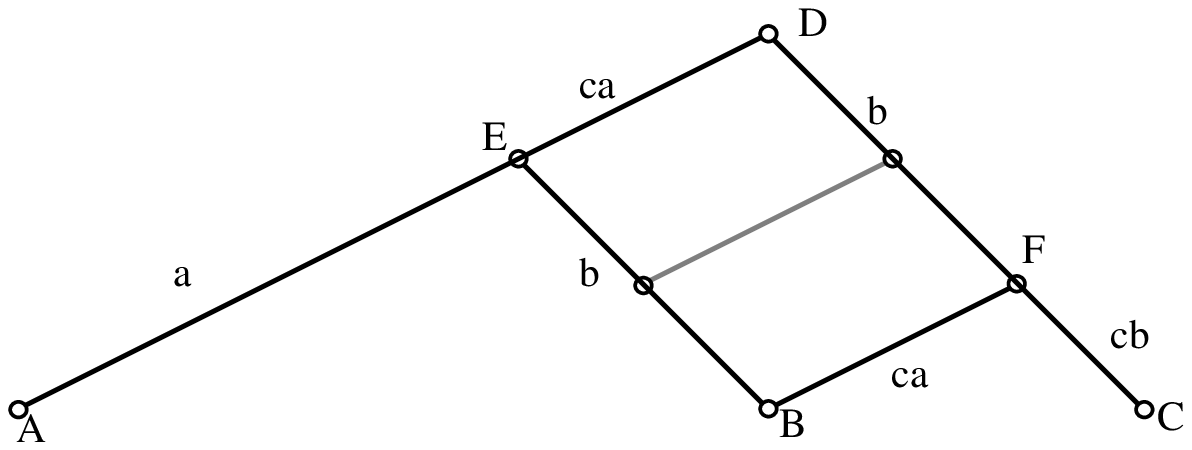 scaled 700}
\capshun {10} {The  Pantograph}
\end{figure}

For any realization $\varphi $ of the pantograph
\begin{equation}\label{eqn:4}
\varphi (C)-\varphi (A)=(1+c)(\varphi (B)-\varphi (A))
\end{equation}
So the pantograph is a quasifunctional
linkage by equation (\ref{eqn:4}).
For example if vertex $A$ is fixed at 0 and
$B$ is the input and $C$ is the output,
it will be quasifunctional for $x\mapsto (1+c)x$.
But if $n>2$, the configuration space is too big for it to be a functional
linkage since any realization can be rotated about
the line through $A$, $B$ and $C$.
Consequently, $\rh{}B$ would not be finite to one, so we could not
get a functional linkage.
To take care of this problem, we will combine the 
pantograph with a linkage from \figsix\ to make it functional.

We divide the construction of a functional linkage for
scalar multiplication into three cases, $\lambda >1$, $0<\lambda <1$, and $\lambda <0$.
The remaining cases $\lambda =0$ or $\lambda =1$ are trivial functions which 
have trivial functional linkages.

If $\lambda >1$ we take $c=\lambda -1$, let $B$ be the input vertex and let
$C$ be the output vertex, and fix $A$ at 0.
We add fixed vertices $A_i$ for $n\ge i\ge 3$, fixed at $ae_i$,
and edges $\lineseg {A_i}E$ of length $\sqrt 2 a$.
So we have a \figsix\ sublinkage $\ul'$ with vertices $A$, $A_i$, $E$, and $B$.
Note that for any $\varphi \in \conf(\ul)$ we have $\varphi (C)=\lambda \varphi (B)$.
So $\ul$ is \qf.
Figure \ref{fig:6} shows this linkage for $n=4$ and $\lambda =1+c$.

\begin{figure}
\Myepsf{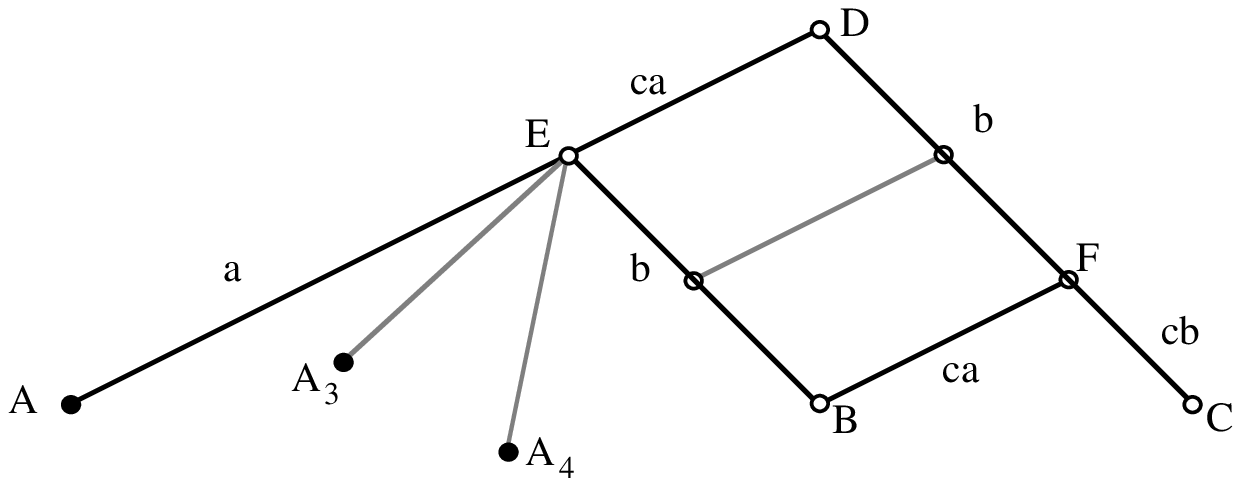 scaled 700}
\capshun 6 {A  functional linkage for scalar multiplication in ${\mathbb R}^4$}
\end{figure}

By Lemma \ref{lem:6.6} we know
that $\rh{}{\ul'}\colon \conf(\ul)\to \conf(\ul')$ is an isomorphism
since the positions of $D$, $F$, and $C$ are polynomial functions
of $A$, $E$, and $B$.
By the discussion of the linkage in \figsix,
we know that if $a$ and $b$ are chosen appropriately,
then $\rh{}B$ double covers some disc $|z-z_0|\le r$.
Hence $\ul$ is functional with \rd\ $\{|z-z_0|\le r\}$.

To get a strong functional linkage, we add two fixed vertices and
tether $E$ and $B$ to them with appropriate length cables
as in \figsix.
By the discussion in \figsix,
we know that $\rh{}B$ singly covers some disc $|z-z_0|\le r$.
Consequently we get a strong functional
linkage with \rd\ $\{|z-z_0|\le r\}$.

If $0<\lambda <1$, we take $c=1/\lambda -1$, let $C$ be the input vertex,
let $B$ be the output vertex, and fix $A$ at 0.
We add fixed vertices $A_i$ at $(a/\lambda )e_i$ 
and edges $A_iD$ of length $\sqrt 2 a/\lambda $ for $3\le i\le n$.
By considering the sublinkage with vertices $A$, $A_i$, $D$, and $C$,
we see as above that with appropriate choices of $a$ and $b$,
the linkage will be functional for $z\mapsto \lambda z$ with \rd\
$|z-z_0|\le r$.
To get a strongly functional linkage, tether $D$ and $C$
appropriately to fixed vertices $A_1$ and $A_2$, as in \figsix.
For the next section it will be useful to point out how $z_0$ can
be chosen.
Looking back at the analysis of \figsix,
we see that we can if we wish pick $z_0$ to be any point
in the circle of radius $a/\lambda $ in the $x_1x_2$ plane.

If $\lambda <0$, we take $c=-\lambda $, let $A$ be the input vertex, 
$C$ be the output vertex, and fix $B$ at 0.
We add fixed vertices $B_i$ at $be_i$ and edges $B_iE$ of length $\sqrt 2 b$.
Letting $\ul'$ be the sublinkage with vertices $B$, $B_i$, $E$, and $A$,
we see as above that $\ul$ is functional for $z\mapsto \lambda z$ with \rd\
$|z-z_0|\le r$.
To get a strongly functional linkage, tether $E$ and $A$
appropriately.

\subsection{A Functional Linkage for the Average}

Now let us find a functional linkage $\ul$ for the average.
Again $\ul$ will be based on the pantograph of Figure \ref{fig:10}.
The input vertices will be $A$ and $C$.
The output vertex will be $B$.  
We let $c=1$ and choose $a=b>r$.
Note $\ul$ is \qf\ for $(z,w)\mapsto (z+w)/2$.
However, this $\ul$ could not be functional because
$D$ is free to rotate around the line through $A$, $B$, and $C$;
hence $(\rh{}C,\rh{}A)$ could not be finite to one.
So we modify the pantograph as follows.

Let us start with the linkage $\ul'''$ found in the previous section 
which is (strongly) functional for the map $z\mapsto z/2$
with restricted domain $U'''=\{ |z-2ae_1|\le 2r \}$.
In particular, $\ul'''$ is a pantograph together with a few more fixed
vertices $A_\ell ,\ldots ,A_n$,  where $\ell =1$ in the \rl\ case
and $\ell =3$ in the classical linkage case.
The fixed vertices are $A$ fixed at $0$, and $A_i$ fixed at some $z_i$.
The input vertex is $C$ and the output is $B$.
Let $\ul''$ be obtained from $\ul'''$ by unfixing all the fixed vertices of $\ul'''$.

Using Lemma \ref{lem:6.9} and the functional linkage found above
for translation, we may find a (strong) functional linkage $\ul'$
for the function $z\mapsto (z+z_\ell ,\ldots ,z+z_n)$
with \rd\ $U'=\{ |z+ae_1|\le r \}$.
Let the input vertex of $\ul'$ be $v$ and let the output vertices
of $\ul'$ be $w_\ell ,\ldots ,w_n$.
We form a linkage $\ul$ by taking the disjoint union of $\ul'$ and $\ul''$,
identifying $v$ with $A$,
and identifying $w_i$ with $A_i$.
Figure \ref{fig:7} shows the result for a classical linkage in ${\mathbb R}^3$.

\begin{figure}
\Myepsf{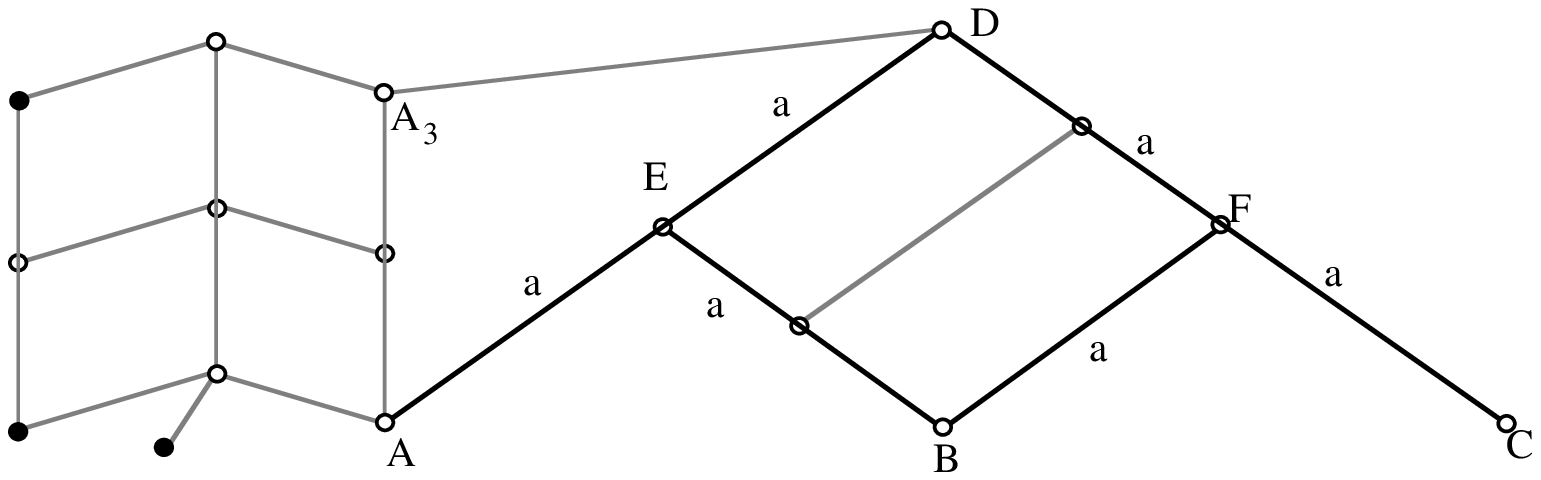 scaled 700}
\capshun 7 {A functional linkage for the average in ${\mathbb  R}^3$}
\end{figure}

By Lemma \ref{lem:6.1} we know that
\begin{eqnarray*}
\conf(\ul)&=&\{ (\varphi ',\varphi '')\in \conf(\ul')\times \conf(\ul'') \mid \varphi '(v)=\varphi ''(A) \mand \varphi '(w_i)=\varphi ''(A_i) \}\\
&=&\{ (\varphi ',\varphi '')\in \conf(\ul')\times \conf(\ul'') \mid \varphi '(v)=\varphi ''(A) \mand \varphi '(v)+z_i=\varphi ''(A_i) \}
\end{eqnarray*}
Thus we have an isomorphism 
$$\beta \colon \conf(\ul')\times \conf(\ul''')\to \conf(\ul)$$
given by $\beta (\varphi ',\varphi ''')=(\varphi ',\varphi '''+\varphi '(v))$  where $\varphi '''+\varphi '(v)$ is the translate
of $\varphi '''$ by $\varphi '(v)$.
Now $$(\rh{}A,\rh{}C)\compose \beta (\varphi ',\varphi ''')=(\varphi '(v),\varphi '''(C)+\varphi '(v))$$
By (strong) functionality we have 
finite sets $F'$ and $F'''$ and analytic isomorphisms
\begin{eqnarray*}
\sigma  '\colon U'\times F' &\to & {\rh{'}v}^{-1}(U')\\
\sigma  '''\colon U'''\times F''' & \to & {\rh{'''}C}^{-1}(U''')
\end{eqnarray*}
so that $\rh{'}v\sigma  '(x,f)=x$ and $\rh{'''}C\sigma  '''(x,f)=x$, and so that
$F'$ and $F'''$ are singletons in the \rl\ case.

In the \rl\ case, we also tether $A$ to $-ae_1$ 
and tether $C$ to $ae_1$ with cables of length $r$.
This will insure that the domain of $\ul$ is $U$ given in (\ref{eqn:7}) below.

So in any case, if 
\begin{equation}\label{eqn:7}
U=\{\,(z,w)\in {\mathbb R}^n\times {\mathbb R}^n \mid 
r\ge |z-ae_1| \mand r\ge |w+ae_1|\,\}
\end{equation}
we have an analytic isomorphism
$$\sigma  \colon U\times F' \times F''' \to   {(\rh{}A,\rh{}C)}^{-1}(U)$$
given by
$$\sigma  (z,w,f',f''')=\beta (\sigma  '(w,f'),\sigma  '''(z-w,f'''))$$
In particular, $(\rh{}C,\rh{}A)$ is an analytically trivial
cover of $U$,
and is an analytic isomorphism in the \rl\ case.
So $\ul$ is (strongly) functional for the average.

As it happens, this is the only place where a \rl\
is constructed with a flexible edge between two nonfixed vertices.
If one wished, one could change the construction slightly
to avoid this, by eliminating $A_1$ and $A_2$,
but tethering $D$ to some fixed vertex instead.
One could then strengthen Theorem \ref{thm:6} to conclude that
in addition, each flexible edge of $\ul'$ is connected to a fixed vertex.

\subsection{A Functional Linkage for Inversion in a Line}\label{sec:4.9}

Let $L$ be a line through 0, $z_0\in L-0$ and $K\subset L-0$ compact.
We will now construct a functional linkage $\ul$ for the function
$f\colon L\to L$
with \rd\ $K$,
where $f(sz_0)= (1/s)z_0$.
In the \rl\ case we will make $\ul$ strongly functional.
Note that the input vertex $v$ will be constrained to lie on the line $L$.
Here we will restrict to the case $n\ge 3$.
The $n=2$ case requires a different construction,
but was shown in \cite{KM} and \cite{K1} so we will not repeat it here.

It suffices to show this for $L$ the $x_1$ axis, and $z_0=e_1$,
 since any other
line and $z_0$ may be obtained by rotation of $\ul$,
and rescaling all side lengths.
Pick $0<c<b$ to be determined later.
Let $a=\sqrt{1+b^2}$ and $d=\sqrt{1+c^2}$.

For $i=2,\ldots ,n$, let $S_i$ be the circle 
$$S_i=\{ x\in {\mathbb  R}^n \mid x_1^2+x_i^2=a^2, \mand x_j=0
\mfor j\neq 1,i\,\}$$
Using Lemma \ref{lem:6.11}, find linkages $\ul_{ij}$ with vertices $v_{ij}$ so that
$\sconf(\ul_{ij},\{v_{ij}\})=S_i$, 
and $\rh{_{ij}}{v_{ij}}$ is an isomorphism
$i=2,\ldots ,k$ and  $j=0,1$.
Form a linkage $\ul$ by taking the disjoint union of the linkages
$\ul_{ij}$, adding vertices $v$ and $w$, and
putting in rigidified squares with vertices $v$, $v_{i0}$, $w$, $v_{i1}$,
and side lengths $b$.
In the \rl\ case put in a cable between $v$ and $w$ with length $2c$.
We also tether each $v_{i0}$ to the point $ae_i$ with a cable of length
$\sqrt 2 a$.
In the classical linkage case, 
 use Lemma \ref{lem:6.11} to find a linkage $\ul_1$ and vertex $v_1$ so that
 $$\sconf(\ul_1,\{v_1\})=\{ x\in {\mathbb  R}^n \mid x_1^2+x_2^2=d^2, \mand x_j=0
\mfor j\ge 3\,\}$$
and $\rh{_1}{v_1}$ is an isomorphism to this circle.
Add $\ul_1$ to $\ul$ and put in edges 
$\lineseg v{v_1}$ and $\lineseg w{v_1}$ of length $c$.

\begin{figure}
\Myepsf{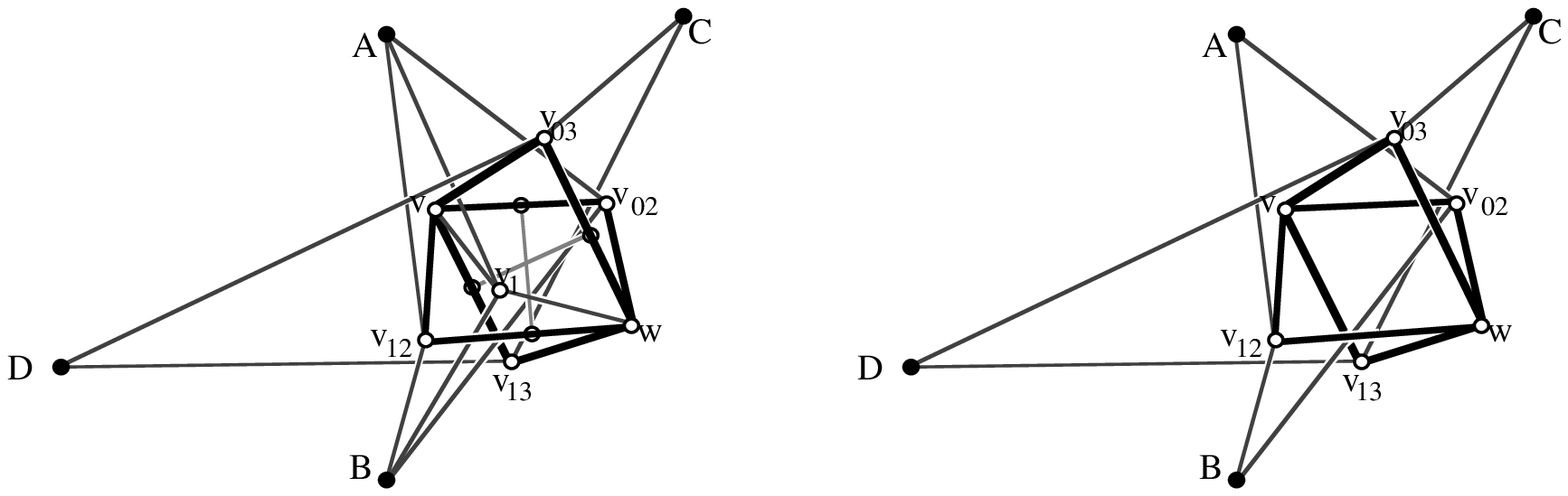 scaled 700}
\capshun 8 {A functional linkage for inversion in a line in ${\mathbb  R}^3$}
\end{figure}

Figure \ref{fig:8} shows a functional classical linkage for 
inversion in the $x_1$ axis of ${\mathbb  R}^3$,
viewed from a point $(e,e,e)$.
At the left is the full linkage which is somewhat complicated
to decipher, the right hand linkage eliminates the rigidifying edges
and the vertex $v_1$ to give the essentials.
Vertices $A$ and $B$ are fixed at points on the $x_3$-axis and
$C$ and $D$ are fixed at points on the $x_2$-axis.
Thus $v_{02}$, $v_{12}$, and $v_1$ are restricted to circles in the
$x_1x_2$ plane, and $v_{03}$ and $v_{13}$ are restricted to a circle
in the $x_1x_3$ plane.
The thick edges form the important part of the linkage,
two rigidified squares.

Let us show that $\ul$ is \qf\ for inversion.
Pick any $\varphi \in \conf(\ul)$ and let $\varphi (v)=x$, $\varphi (w)=y$ 
and $\varphi (v_{ij})=x_{ij}$.
If $x_{i0}= x_{i1}$ then $|x-y|=2b>2c$.
So in the \rl\ case the cable between $v$ and $w$ guarantees that 
 $x_{i0}\neq x_{i1}$.
In the classical case, $|x-y|\le |x-\varphi (v_1)|+|\varphi (v_1)-y|=2c$,  so again
 $x_{i0}\neq x_{i1}$.

Because we rigidified the quadrilateral $x x_{i0} y x_{i1}$, we know that
$(x+y)/2=(x_{i0} + x_{i1})/2$ since these are both midpoints of the 
intersecting diagonals.
So $x+y$ is in the $x_1x_i$ plane for all $i\ge 2$.
Since $n\ge 3$, we then conclude that $x+y$ is on the $x_1$ axis.
So 
$(x+y)/2=\alpha e_1$ for some $\alpha $.
Since $\alpha e_1$ is at the midpoint of a chord of the circle $S_i$,
we must have $x_{ij}=\alpha e_1+\epsilon _i(-1)^j\beta e_i$ where 
$\beta =\sqrt{a^2-\alpha ^2}>0$
and $\epsilon _i=\pm 1$.
Note that in the \rl\ case that 
 the cable tethering
$v_{i0}$ to $ae_i$ means that 
\begin{equation}\label{eqn:10}
2a^2\ge |\alpha e_1+\epsilon _i\beta e_i-ae_i|^2=
\alpha ^2+\beta ^2-2\epsilon _i\beta a+a^2=2a^2-2\epsilon _i\beta a
\end{equation}
 so $\epsilon _i\beta a\ge 0$ and hence  $\epsilon _i=1$.

Now
\begin{equation}\label{eqn:1}
b^2=|x-x_{ij}|^2=|x|^2+|x_{ij}|^2-2 x\cdot x_{ij}=|x|^2+a^2-2x\cdot x_{ij}
\end{equation}
Consequently $x\cdot (x_{i0}-x_{i1})=0$ 
so $x \cdot e_i =0$ 
for all $i\ge 2$.
So $x=\gamma  e_1$ for some $\gamma $.
Note from equation (\ref{eqn:1}) that $b^2=\gamma ^2+a^2-2\gamma \alpha $,
so $\gamma ^2-2\gamma \alpha +1=0$. 
Hence
\begin{eqnarray*}
\gamma &=&\alpha \pm \sqrt{\alpha ^2-1}\\
y&=&2(x+y)/2-x=(2\alpha -\gamma )e_1=(1/\gamma  ) e_1
\end{eqnarray*}
So  we see that $\ul$ is
\qf\ for $se_1\mapsto (1/s)e_1$.

Let us now see what the \rd\ can be.
Since $|x-y|\le 2c$ we must have $|\gamma -1/\gamma |\le 2c$
so solving we find we must have $\gamma \in A$ where 
 $$A=\{ t\in {\mathbb  R} \mid d-c\le |t|\le d+c\}$$
 By choosing $c$ large enough, we may ensure that
 $K\subset Ae_1$.

We will show that  
 the image of 	$\rh{}v$ is $Ae_1$.
 Moreover in the \rl\ case $\rh{}v$ is an analytic isomorphism to $Ae_1$,
 and in the classical case $\rh{}v$ restricts to an analytically trivial  cover
 of  $A'e_1$ where $A'$ is the interior of $A$.
 
So we need to solve for the positions of the vertices in terms of $\gamma $.
Pick any $\gamma \in A$.
 We have already seen that if $\varphi \in \conf(\ul)$ and $\varphi (v)=\gamma e_1$,
 then $\varphi (w)=(1/\gamma )e_1$.
 Moreover, if $\alpha =(\gamma +1/\gamma )/2$ and $\beta =\sqrt{a^2-\alpha ^2}$ then 
\begin{equation}\label{eqn:9}
\varphi (v_{ij})=\alpha e_1+\epsilon _i(-1)^j\beta e_i
\end{equation}
 where  $\epsilon _i=\pm 1$,
and $\epsilon _i=1$ in the \rl\ case.
In the \rl\ case, the cable between $v$ and $w$ means we must have
$|\gamma -1/\gamma |\le 2c$ which is true for all $\gamma \in A$.
Also the cables between $ae_i$ and $v_{i0}$ require that
$|\varphi (v_{i0})-ae_i|\le \sqrt 2 a$ which follows from (\ref{eqn:10}) since 
$\epsilon _i=1$.
 
 So in the \rl\ case we have seen that $\rh{}v$ is an analytic isomorphism 
from $\conf(\ul)$ to $Ae_1$, hence $\ul$ is strongly functional for $f$
with domain $Ae_1\supset K$.

In the classical case,  it only remains to solve for $\varphi (v_1)$.
We know that $\varphi (v_1)$ must be in the circle of radius $d$ in the $x_1x_2$ plane
which means that $\varphi (v_1)=se_1+te_2$ for some $s,t$ with $s^2+t^2=d^2$.
So we must have
\begin{eqnarray*}
0=|\varphi (v_1)-\varphi (v)|^2-c^2&=& (s-\gamma )^2+t^2-c^2\\
&=&d^2-c^2-2s\gamma +\gamma ^2=1-2s\gamma +\gamma ^2\\
0=|\varphi (v_1)-\varphi (w)|^2-c^2&=& (s-1/\gamma )^2+t^2-c^2\\
&=&d^2-c^2-2s/\gamma +1/\gamma ^2=1-2s/\gamma +1/\gamma ^2
\end{eqnarray*}
From either of these equations we then solve for $s$ and obtain
$s=\alpha $.
Consequently, 
\begin{equation}\label{eqn:8}
\varphi (v_1)=\alpha e_1+\epsilon _1\sqrt{d^2-\alpha ^2}\,e_2
\end{equation}
 where $\epsilon _1=\pm 1$.

So in the end we see we have a map
$\sigma  \colon A\times \{-1,1\}^n\to \conf(\ul)$ 
where the $\{-1,1\}^n$ chooses the signs $\epsilon _i$.
Moreover, looking at equations (\ref{eqn:8}) and 
(\ref{eqn:9}) and recalling that $a>d$, we see that $\sigma $ restricts to
an analytic isomorphism  wherever we have
$0<d^2-\alpha ^2=c^2-(\gamma -1/\gamma )^2$,
i.e., on $A'\times \{-1,1\}^n$.
So $\ul$ is functional for $f$ with \rd\ $A'e_1$.

This construction of inversion in a line is the one place in this
paper where we have used $n\ge 3$.
If $n=2$, a different construction is needed,
see \cite{K1} or \cite{KM}.
The construction above with $n=2$ would actually give you 
 inversion through the circle, since the domain would not be restricted
to a single line.

\subsection{A Functional Linkage for Orthogonal Projection to a Line}

Let $g\colon {\mathbb  R}^n\to {\mathbb  R}^n$ be orthogonal projection to a line $L$.
Let $K\subset {\mathbb  R}^n$ be compact.
We will now construct a functional linkage $\ul$ for $g$ with \rd\ $K$.
As a bonus, we will at the same time construct a functional
linkage for reflection $f\colon {\mathbb  R}^n\to {\mathbb  R}^n$
about the line $L$, although we do not use this fact.

Before finding this linkage, we prove the following Lemma:

\begin{lem}\label{lem:6.12}
Given any compact line segment $T\subset {\mathbb  R}^n$, there is a linkage $\ul$ and a
$v\in \vrtcs(\ul)$ so that $T=\sconf(\ul,\{v\})$.
If we insist that $\ul$ be a classical linkage then
we may ensure that $\rh{}v$ restricts to an analytically trivial cover of
the interior $T'$ of $T$.
If $\ul$ is allowed to be a  \rl, we may make 
$\rh{}v\colon \conf(\ul)\to T$ be an analytic isomorphism.
\end{lem}

\begin{proof}
It suffices to prove this for only one line segment $T$,
since any other line segment may be obtained from $T$ by
translation, rotation, and rescaling.
Let $\ul_0$ be the (strong) functional linkage constructed 
in section \ref{sec:4.9} for
inversion in the $x_1$ axis.
Take $c=3/4$ and $d=5/4$ so the domain is $U_0=\{ se_1\mid 1/2\le |s|\le 2\}$.
Let $\ul_1$ be obtained from $\ul_0$ by translating by $2e_1$,
so the domain of $\ul_1$ is $U_1=\{ se_1\mid 1/2\le |s-2|\le 2\}$.
Let $v_i$ be the input vertex of $\ul_i$.
Let $\ul$ be obtained from the disjoint union of $\ul_0$ and $\ul_1$
by gluing the inputs $v_0$ and $v_1$ together.
Then $\conf(\ul)$ is the fiber product of the maps $\rh{_0}{v_0}$
and $\rh{_1}{v_1}$.
Consequently, if $v=v_0=v_1$ is the glued vertex, we have
$$\rh{}v(\conf(\ul))=U_0\cap U_1
=\{ se_1\mid 1/2\le s\le 3/2\}=T$$

Let $U'_i$ and $T'$ be the interiors of $U_i$ and $T$ in the $x_1$ axis.
We have a finite set $F$ and analytic isomorphisms
$\sigma  _i\colon U_i'\times F\to  \rh{_i}{v_i}^{-1}(U_i')$ so that $\rh{_i}{v_i}(\sigma  _i(u,f))=u$.
Looking at $\conf(\ul)$ as a fiber product
$$\conf(\ul)=\{ (\varphi _0,\varphi _1)\in \conf(\ul_0)\times \conf(\ul_1) 
\mid \varphi _0(v_0)=\varphi _1(v_1) \}$$
we have an analytic isomorphism
$\sigma  \colon T'\times F\times F\to \rh{}{v}^{-1}(T')$
given by $\sigma  (u,f_0,f_1)=(\sigma  _0(u,f_0), \sigma  _1(u,f_1))$.

In the \rl\ case, we may take each $\rh{_i}{v_i}$ to be an analytic
isomorphism and hence $\rh{}{v}\colon \conf(\ul)\to T$ is an 
analytic isomorphism also.
\end{proof}

We now proceed to find a functional linkage for
projection $g$ to a line $L$.
After translation and rotation, we may as well assume that
$L$ is the $x_1$ axis, and so  $g(x)=x_1e_1$.
Choose $r$ so that $K\subset \{|z|<r/2\}$.
Let $U\subset {\mathbb  R}^n$ be the set of points of distance $\le r$ from $L$.
Define $\psi _i\colon U\to L$, $i=0,1$ by 
$$\psi _i(x)=g(x)+(-1)^i\sqrt{r^2-|x-g(x)|^2}\,e_1$$
Note that $\psi _i(x)$ are the two points on $L$ with distance $r$
from $x$.
For $x\in K$ we must have $|g(x)|\le |x|<r/2$ and
$\sqrt{r^2-|x-g(x)|^2}> \sqrt{r^2-(r/2)^2}=\sqrt 3 r/2$.
So if $L_0=\{ te_1 \mid t>0\}$ and $L_1=\{te_1 \mid t<0\}$
we must have $\psi _i(K)\subset L_i$.

Choose  closed line segments $T_i$ in $L_i$ so that
$\psi _i(K)$ is contained in the interior $T'_i$ of $T_i$, $i=0,1$.
By Lemma \ref{lem:6.12}, we may choose linkages $\ul_0$ and $\ul_1$ with
vertices $v_i$ so that $\sconf(\ul_i,\{v_i\})=T_i$, and
$\rh{_i}{v_i}\colon \conf(\ul_i)\to L$ is an analytically trivial covering
of $T'_i$, and in fact in the \rl\ case it is an analytic isomorphism
to $T_i$.

Form $\ul$ by taking the disjoint union of $\ul_0$ and $\ul_1$,
adding two vertices $v$ and $w$,
and forming a rigidified square $v_1vv_0w$ with side length $r$.
To get a functional linkage for $g$,
we will be a bit specific about this rigidification.
We will add three vertices, $v_2$, $v_3$, $v_4$.
We will place $v_4$ at the midpoint of the edge $\lineseg {v_1}v$,
place $v_3$ at the midpoint of the edge $\lineseg {v_0}w$,
and place $v_2$ at the midpoint of the edge $\lineseg {v_3}{v_4}$.
Thus, in any realization, $v_2$ will be at the exact center
of the parallelogram $v_1vv_0w$.

\begin{figure}
\Myepsf{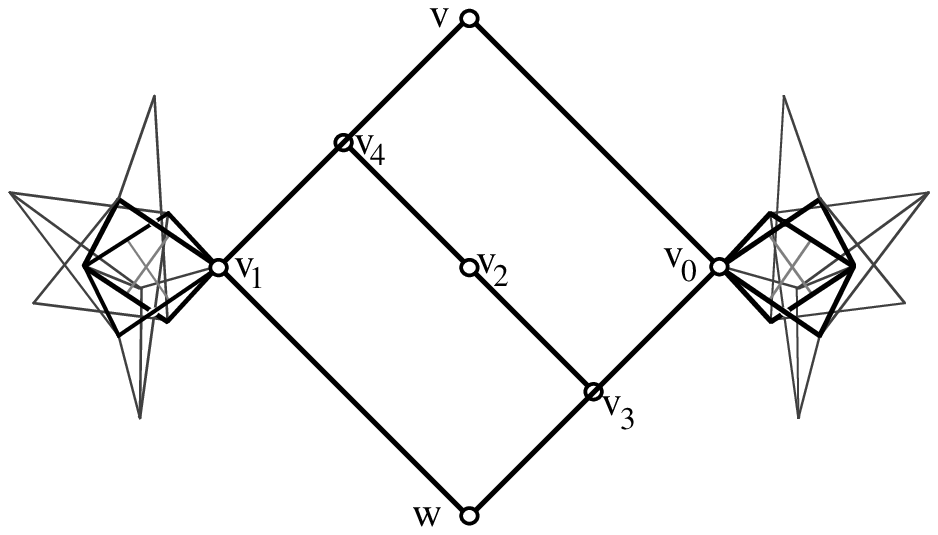 scaled 700}
\capshun 9 {A functional linkage for  projection to a line in ${\mathbb  R}^3$}
\end{figure}

For a functional linkage for projection $g$,
set $v$ to be the input vertex and $v_2$ to be the output vertex.
For a functional linkage for reflection $f$,
set $v$ to be the input vertex and $w$ to be the output vertex.

Pick any $\varphi \in \conf(\ul)$ and let
 $x=\varphi (v)$.
Since $\varphi |_{\vrtcs(\ul_i)}\in \conf(\ul_i)$, we know that $\varphi (v_i)\in T_i\subset L_i$.
Since $T_0$ and $T_1$ are disjoint, we know that $\varphi (v_0)\neq \varphi (v_1)$.
Consequently $\varphi (v_i)=\psi _i(x)$.
Note that 
\begin{eqnarray*}
\varphi (w)&=&\varphi (v_0)+\varphi (v_1)-\varphi (v)
=\psi _0(x)+\psi _1(x)-x=2g(x)-x=f(x)\\
\varphi (v_2)&=&(\varphi (v)+\varphi (w))/2=(x+2g(x)-x)/2=g(x)
\end{eqnarray*}
so $\ul$ is quasifunctional for $f$ and $g$.

Let $\sigma  _i\colon T'_i\times F_i\to \rh{_i}{v_i}^{-1}(T'_i)$ be analytic isomorphisms
with $\rh{_i}{v_i}(\sigma  _i(z,f_i))=z$.
Then we get an analytic isomorphism $\sigma  \colon K\times F_0\times F_1\to \rh{}v^{-1}(K)$
given by 
\begin{eqnarray*}
\sigma  (x,f_0,f_1)|_{\ul_i}&=&\sigma _i(\psi _i(x),f_i)\\
\sigma (x,f_0,f_1)(v)&=&x\\
\sigma (x,f_0,f_1)(w)&=&f(x)\\
\sigma (x,f_0,f_1)(v_2)&=&g(x)\\
\sigma (x,f_0,f_1)(v_3)&=&(f(x)+\psi _0(x))/2\\
\sigma (x,f_0,f_1)(v_4)&=&(x+\psi _1(x))/2
\end{eqnarray*}
So $\ul$ is functional for $f$ and $g$ with \rd\ $K$,
and in the \rl\ case it is strongly functional since the $F_i$
have just one point each and $\sigma $ may be extended to an analytic isomorphism
$\sigma \colon (\psi _0^{-1}(T_0)\cap \psi _1^{-1}(T_1))\times F_0\times F_1\to \conf(\ul)$.

The astute reader will notice that the linkage shown in Figure \ref{fig:9}
is simpler than that constructed in the text (which would 
be hopelessly cluttered with
two copies of the linkage from Figure \ref{fig:8} on each side).
The linkage shown will still work, the only difference is that
each $T_i$ will be a double interval.

\section{Proofs of Theorems}

Now that Theorem \ref{thm:6} is proven, we can prove the theorems 
stated in the first section.  
First we prove a special case of Theorem \ref{thm:3}.

\begin{prop}\label{prop:6.14}
Suppose $Z\subset ({\mathbb  R}^n)^k$ is  compact, $n\ge 2$.
The following are equivalent:
\begin{enumerate}
\item  There is a classical linkage $\ul$ and  a $W\subset \vrtcs(\ul)$ so that
$\sconf(\ul,W)=Z$.
\item  $Z$ is a semialgebraic set.
\end{enumerate}
\end{prop}

\begin{proof}
The implication 1 implies 2 follows 
from the Tarski-Seidenberg Theorem, \cite{S},
because $\conf(\ul)$ is an algebraic set
and $\sconf(\ul,W)$ is its image under projection.

Now let us show that 2 implies 1.
By Lemma 3.1 of \cite{K2} that there is a polynomial map
$q\colon Y\to ({\mathbb  R}^n)^k$ from some compact real algebraic set $Y$
so that $Z=q(Y)$.
By taking the graph of $q$, we may as well assume that 
$Y\subset ({\mathbb  R}^n)^{k+m}$ and 
$q$ is projection to the first $nk$ coordinates.
Pick a polynomial $p\colon ({\mathbb  R}^n)^{k+m}\to {\mathbb  R}$ so that $Y=p^{-1}(0)$.
By Theorem \ref{thm:6}, there is a functional classical linkage $\ul'$ for the map 
$x\mapsto p(x)e_1$, with \rd\ $Y$
and with distinct input vertices.
Let $\ul$ be obtained from $\ul'$ by fixing the output
vertex at 0.
Let $W$ be the first $k$ input vertices of $\ul'$
and let $U$ be all the input vertices of $\ul'$.
Then by Lemma \ref{lem:6.4}, $\sconf(\ul,U)=Y$, so $\sconf(\ul,W)=Z$.
\end{proof}

We are now able to prove  Theorems \ref{thm:4x} and \ref{thm:4}.

\begin{proof} (of  Theorem \ref{thm:4x})
Note that part 3 implies part 2 since any $m$ points of ${\mathbb  R}^n$
are contained in some affine subspace with dimension $m-1$.
Note also that $\conf(\ul)=\sconf(\ul,\vrtcs(\ul))$
so it suffices to prove the $\sconf(\ul,W)$ results only.

To see part 3, let $b=\dim T$.
We may pick $\beta \in \Eu(n)$ so that $\beta (T)={\mathbb  R}^{b}\times 0\subset {\mathbb  R}^{b}\times {\mathbb  R}^{n-b}$.
Then if  $\gamma \in O(n-b)$ we have
$\beta ^{-1}\gamma \beta \ul=\ul$.
So $\sconf(\ul,W)$ is invariant under the conjugate of $O(n-b)$ by $\beta $.
We know $\conf(\ul)$ is a quasialgebraic set, so it is closed.
But it is also
contained in a
ball of radius $d$ around the image of a fixed vertex, where
$d$ is the sum of the lengths of all edges of $\ul$.
(We are using connectedness here.)
So we see that $\conf(\ul)$ is compact.
So then $\sconf(\ul,W)$ is compact, since it is the image of
the compact $\conf(\ul)$ under projection.

To see part 1, note that by Lemma \ref{lem:6.3},
 $\sconf(\ul,W)$ is invariant under the action of $\Eu(n)$.
If $W$ is empty, then $\sconf(\ul,W)$ is a single point which is compact,
so we may assume that $W$ has $k>0$ vertices.
Recall $\Tran(n)\subset \Eu(n)$ is the subgroup of translations.
We may identify the quotient $\sconf(\ul,W)/\Tran(n)$ with
$$Z_0=\{ (z_1,\ldots ,z_k)\in \sconf(\ul,W) \mid z_k=0 \}$$
This is compact by part 2, since it is $\sconf(\ul',W)$ for the
linkage $\ul'$ formed from $\ul$ by fixing the $k$-th vertex
in $W$ to $0$.
Consequently, the quotient $\sconf(\ul,W)/\Eu(n)=Z_0/O(n)$ is compact.
\end{proof}

\begin{proof} (of  Theorem \ref{thm:4})
If $Z$ is empty, we may easily prove this by choosing any $\ul$
which includes a triangle which violates the triangular inequality.
So we may assume $Z$ is nonempty.

Let us first prove parts 2 and 3.
First note that in part 2, by replacing $Z$ by $\beta (Z)$ for some $\beta \in \Eu(n)$,
 we may assume that $G=O(m)$, acting on the last $m$ coordinates of ${\mathbb  R}^n$.
By Proposition \ref{prop:6.14},
  we may find a linkage $\ul'$ and a $W\subset \vrtcs(\ul')$ so that
$\sconf(\ul',W)=Z$.
Throw away all connected components of $\ul'$ which do not contain
any vertices in $W$ or any fixed vertices.
By Lemma \ref{lem:6.2}, doing so does not change $\sconf(\ul',W)$.
By adding some isolated fixed vertices to $\ul'$
if necessary, we may assume that there is a vertex fixed at 0,
a vertex fixed at each $e_i$, $i=1,\ldots ,n$ and a vertex fixed at 
$\sum _{i=1}^n e_i$.
Adding an isolated fixed vertex to $\ul'$ does not change
$\sconf(\ul',W)$.

Let the fixed vertices of $\ul'$ be $\{v_0,\ldots ,v_b\}$
where $v_i$ is fixed to the point $z_i$.
We may suppose $z_0=0$, $z_{i}=e_i$,  $i=1,\ldots n$,
 and $z_{n+1}=\sum _{i=1}^n e_i$.
For each pair $i,j$ with $z_i\neq z_j$ put in an edge $v_iv_j$ of length $|z_i-z_j|$,
if it is not already there.
This will not change $\sconf(\ul',W)$.
Note we did not attempt to add any zero length edges,
which would not be allowed.

Let $\ul''$ be obtained from $\ul'$ by only fixing the vertices
$v_i$ for $i\le n$.
We claim that $\sconf(\ul',W)=\sconf(\ul'',W)$.
One inclusion $\sconf(\ul',W)\subset \sconf(\ul'',W)$ is trivial.
So let us see the other inclusion.
Pick any $\varphi \in \conf(\ul'')$. 
We claim that in fact $\varphi (v_i)=z_i$ for all $i$.
To see this, note first that two different points in ${\mathbb  R}^n$ can not
have the same distances from $n+1$ points in general position,
see Lemma \ref{lem:d}.
Here you can interpret general position to mean that
their convex hull has nonempty interior.
Consequently $\varphi '(v_{n+1})=z_{n+1}$ since the $n+1$ edges 
$\lineseg {v_i}{v_{n+1}}$, $i\le n$
have lengths $|z_{n+1}-z_i|$, so 
$$|\varphi '(v_{n+1})-z_i|=|\varphi '(v_{n+1})-\varphi '(v_i)|=|z_{n+1}-z_i|$$
For any $j>{n+1}$, there are edges in $\ul$ from $v_j$ to at least $n+1$
of the $v_i$, $i\le {n+1}$, and consequently $\varphi '(v_j)=z_j$ since any $n+1$
of the $z_i$, $i\le {n+1}$ are in general position.
Consequently, $\varphi \in \conf(\ul')$.
So $\varphi |_{W}\in \sconf(\ul',W)$, and we have shown that
$\sconf(\ul',W)=\sconf(\ul'',W)$.

We claim that $\ul''$ is connected.
Note $Z$ is compact since its quotient by $O(m)$ is compact.
By Lemma \ref{lem:6.13}, we see that each
connected component of $\ul''$ has at least one fixed
vertex.
But the $n+1$ fixed vertices are connected to each other, hence 
there is only one
connected component, so part 3 is proven.

Let us now prove part 2.
Let $\ul$ be obtained from $\ul''$ by only fixing the vertices $v_i$
for $i\le n-m$,
and not fixing any of the other vertices of $\ul'$.
We claim that $\sconf(\ul'',W)=\sconf(\ul,W)$.
Again, one inclusion $\sconf(\ul'',W)\subset \sconf(\ul,W)$ is trivial.
So let us see the other inclusion.
Pick any $\varphi \in \conf(\ul)$. 
Now if $1\le i\le n$ and $1\le j\le n$, $j\neq i$
we have
$|\varphi (v_i)|=|z_i|=1$, and $|\varphi (v_i)-\varphi (v_j)|=|z_i-z_j|=\sqrt 2$.
Consequently, the vectors $\varphi (v_1),\ldots ,\varphi (v_{n})$
form an orthonormal set.
Since $\varphi (v_i)=e_i$ for $i\le n-m$, we may choose a $\beta \in O(m)$,
acting on the last $m$ coordinates of ${\mathbb  R}^n$,
so that $\beta \varphi (v_i)=z_i$ for all $i\le n$.
For convenience, let $\varphi '=\beta \compose \varphi $.
Note that $\varphi '\in \conf(\ul'')$.
So $\varphi '|_{W}\in \sconf(\ul'',W)=Z$.
By $O(m)$ invariance of $Z$, we know that 
$\beta ^{-1}\compose \varphi '|_{W}\in Z$ also.
But $\beta ^{-1}\compose \varphi '|_{W}=\varphi |_{W}$, so $\varphi |_{W}\in Z$.
So we have shown that $\sconf(\ul,W)\subset \sconf(\ul'',W)$,
and hence $\sconf(\ul,W)=Z$.
So part 2 is proven.

Now let us prove part 1.
Let $\Tran(n)\subset \Eu(n)$ be the subgroup of translations.
Then we may identify the quotient $Z/\Tran(n)$ with
$$Z_0=\{ (z_1,\ldots ,z_k)\in Z \mid z_k=0 \}$$
Note $Z_0$ is invariant under the diagonal action of $O(n)$
and $Z_0/O(n)=Z/\Eu(n)$.
So $Z_0/O(n)$ is compact and hence $Z_0$ is compact.
Also $Z_0$ is semialgebraic since it is the intersection of semialgebraic sets.
Define $Z_0'\subset ({\mathbb R}^n)^{k-1}$ by  $Z_0=Z_0'\times 0$.
Note $Z_0'$ is compact and semialgebraic since it is a projection of $Z_0$.
Also $Z_0'$ is $O(n)$ invariant since $Z_0$ is.
So by part 2, there is a linkage $\ul'$ and a $W'\subset \vrtcs(\ul')$
so that $\sconf(\ul',W')=Z_0'$.
Furthermore, there is only one
fixed vertex $v_0$ of $\ul'$ and that vertex is fixed at the point 0.

Consider the linkage $\ul$ formed from $\ul'$ by unfixing the
vertex $v_0$.
So $\ul$ has no fixed vertices.
Let $W=\{w_1,\ldots ,w_k\}$ be such that $W'=\{w_1,\ldots ,w_{k-1}\}$
and $w_k=v_0$.
Then by Theorem \ref{thm:4x}, 
$\sconf(\ul,W)$ is invariant under the action of $\Eu(n)$,
as is $Z$.
But then $\sconf(\ul,W)=Z$ since they both have the same intersection 
with $({\mathbb  R}^n)^{k-1}\times 0$, (namely $Z_0$),
and so are both the union of $\Tran(n)$ orbits of the
same set.
\end{proof}



We may now prove Theorem \ref{thm:3}.

\begin{proof}
It is immediate that condition 2 implies condition 1.
Now let us see why condition 1 implies condition 3.
By Lemma \ref{lem:6.2} it suffices to consider the case where $\ul$
is connected, but this case follows from Theorem \ref{thm:4x}.
To see condition 3 implies condition 2, note that by Lemma \ref{lem:6.2} it suffices
to consider the case where $X$ is itself virtually compact.
But then this case follows from Theorem \ref{thm:4}.
\end{proof}

Now let us prove Theorems \ref{thm:2} and \ref{thm:1} together.

\begin{proof}
Suppose $\ul$ is a linkage.
We need to show that $\conf(\ul)$ is isomorphic to 
$X\times ({\mathbb R}^n)^k$ for some compact quasialgebraic set $X$,
and $X$ is algebraic if $\ul$ is classical.
By Lemma \ref{lem:6.2}, we may suppose that $\ul$ is connected.
If $\ul$ has any fixed vertices, then $\conf(\ul)$ is  compact
by Theorem \ref{thm:4x}.
So we have reduced to the case where $\ul$ has no fixed vertices.
Form a linkage $\ul'$ from $\ul$ by fixing one of its vertices $v$ to $0$.
By Theorem \ref{thm:4x}, $\conf(\ul')$ is compact.
By Lemma \ref{lem:3.4}, $\conf(\ul)$ is isomorphic to 
$\conf(\ul')\times {\mathbb R}^n$.
So the first parts of the two theorems are proven.

Now let us prove the second part of Theorem \ref{thm:2}.
Let $X\subset ({\mathbb R}^n)^m$ be a compact quasialgebraic set and $k\ge 0$.
Find polynomials $r_i\colon ({\mathbb R}^n)^m\to {\mathbb R}$, $i=1,\ldots ,\ell$ and
$\ell '\le \ell$ so that 
\begin{eqnarray}\label{eqn:2}
X&=&\{\, x\in ({\mathbb R}^n)^m \mid r_i(x)=0 \mfor i\le \ell ' \mand r_i\ge 0 
\mfor \ell '<i\le \ell \,\}
\end{eqnarray}
Let $r\colon ({\mathbb R}^n)^m\to ({\mathbb R}^n)^\ell $ be the map 
$(r_1e_1,r_2e_1,\ldots ,r_\ell e_1)$.
Let $\ul'$ be a strong functional linkage for the map $r$ with \rd\
 $X$.
 Let its input and output maps be $q$ and $p$ respectively.
 Let $v_1,\ldots ,v_{\ell }$ be the output vertices.
 By compactness we may pick a $d>0$ so that $r_i(x)\le 2d$ for all $x\in X$,
 $i=1,\ldots ,\ell $.
 
 Form a linkage $\ul$ as follows.
 Take $\ul'$ and add a vertex $u_i$ for each $\ell '<i\le \ell $.
 For each $i\le \ell '$ fix the vertex $v_i$ to $0$.
 For each $\ell '<i\le \ell $ attach a flexible edge of length $d$
 between $v_i$ and $u_i$, and fix $u_i$ to $de_1$.
 
 By Lemmas \ref{lem:6.4} and \ref{lem:6.5}, $\conf(\ul)$
 is naturally isomorphic to 
 $$Y=\{ \varphi \in \conf(\ul') \mid \varphi (v_i)=0 \mfor i\le \ell ' \mand |\varphi (v_i)-de_1| \le  d \mfor 
 \ell '< i \le \ell \}$$
 but $p(\varphi )=(\varphi (v_1),\ldots ,\varphi (v_\ell ))$ and $p=r\compose q$, 
 so by (\ref{eqn:2}), $Y=q^{-1}(X)$.
But $q$ is an analytic isomorphism, so $\conf(\ul)$ is analytically
isomorphic to $X$.
Since this isomorphism $q\compose \rh{}{\ul'}$ is just given by coordinate projection, we
immediately obtain the map $\alpha $.

So we have done the compact case, $k=0$.
In general if $k>0$, just add $k$ isolated unfixed vertices to $\ul$.
By Lemma \ref{lem:6.2} the resulting linkage
will have configuration space $\conf(\ul)\times ({\mathbb R}^n)^k$.

The second part of Theorem \ref{thm:1} is proven similarly.
We start with a compact algebraic set $X$.
There are no inequalities, so $\ell '=\ell$.
However, since we want the linkage $\ul'$ to be classical,
it may not be strongly functional, just functional.
So we have an analytic isomorphism $\sigma \colon X\times F\to q^{-1}(X)$
for some finite set $F$.
We form $\ul$ by fixing all the output vertices of $\ul'$ to $0$.
Then by Lemma \ref{lem:6.4} we know that
$$\conf(\ul)=p^{-1}(0)=q^{-1}r^{-1}(0)=q^{-1}(X)$$
and the result follows.
\end{proof}

\begin{proof} (of  Theorem \ref{thm:5})
As in the proof of Theorem \ref{thm:4},
we may assume that $Z$ is nonempty.

Let us  prove part 3.
By Theorem \ref{thm:1} there is a classical linkage $\ul'$,
 a finite set $F$, and an analytic 
$\alpha \colon Z\times F\to ({\mathbb R}^n)^m$ so that the map
$(x,f)\mapsto (x,\alpha (x,f))$ is an analytic 
 isomorphism from $Z\times F$ to $\conf(\ul')$.
 Note that each component of $\ul'$ must have a fixed vertex
 by Lemma \ref{lem:6.13} with $W=\vrtcs(\ul)$.
We now proceed as in the proof of Theorem \ref{thm:4}.
By adding some isolated fixed vertices to $\ul'$
if necessary, we may assume that there is a vertex $v_0$ fixed at 0,
a vertex $v_i$ fixed at each $e_i$, $i=1,\ldots ,n$ and a vertex fixed at 
$\sum _{i=1}^n e_i$.
Add edges of the appropriate length between the fixed vertices.
Finally, unfix all fixed vertices except $v_0,v_1,\ldots ,v_n$.
Just as in the proof of Theorem \ref{thm:4},
the resulting classical linkage $\ul$ has $\conf(\ul)=\conf(\ul')$.

Let us now prove part 1.
After rotation and translation, we may assume $G=O(1)$.
So $G$ has just one nontrivial element $\tau $, reflection about the
hyperplane $\{x_n=0\}$.
By part 3, there is a classical linkage $\ul'$ with  $n+1$ fixed vertices
$v_0,\ldots ,v_n$ fixed at 0 and $e_i$
so that $\conf(\ul')$ is analytically isomorphic to $Z\times F$
for some finite $F$.
Let $\ul$ be obtained from $\ul'$ by unfixing the vertex
$v_n$, but adding edges $\lineseg {v_n}{v_0}$ of length 1 
and $\lineseg {v_n}{v_i}$ of length $\sqrt 2$ for $1\le i\le n-1$,
if these are not already present.
Pick any  $\varphi \in \conf(\ul)$.
Then either $\varphi (v_n)=e_n$ or $\varphi (v_n)=-e_n$.
If $\varphi (v_n)=e_n$ then $\varphi \in \conf(\ul')$.
If $\varphi (v_n)=-e_n$ then $\tau \varphi \in \conf(\ul')$.
Thus we have an isomorphism $\conf(\ul')\times G\to \conf(\ul)$
given by $(\varphi ,g)\mapsto g\varphi $.
So part 1 is shown.

Let us now prove part 2.
After translation and rotation, we may as well assume that
$G=O(2)$.
If $\beta \colon Z'\times SO(2)\to Z$ is the map $\beta (z,g)=gz$, 
let $\beta ^{-1}=(\eta ,\gamma )$ for
entire rational functions $\eta \colon Z\to Z'$ and $\gamma \colon Z\to SO(2)$.
Let 
$$Z''=\{(z,w)\in Z\times {\mathbb R}^n \mid w=\gamma (z)e_{n-1} \}$$  
Note that $Z'\times e_{n-1}\subset Z''$ and is invariant under $O(1)$.
So by part 1 
 we have a connected classical linkage $\ul'$ with only $n$ fixed vertices
 $v_0,\ldots ,v_{n-1}$, fixed at $0$ and $e_i$, and so that
 there is a finite $F$ and 
 an analytic function
$\alpha '\colon Z'\times F\to ({\mathbb  R}^n)^m$ so that 
$$\conf(\ul')=\{\, (x,e_{n-1},\alpha '(x,f)) 
\mid x\in Z' \mand\ f\in F\,\}$$
and so that the map $(x,f)\mapsto (x,e_{n-1},\alpha '(x,f))$ is an analytic isomorphism.

We obtain $\ul$ by just unfixing the vertex $v_{n-1}$, 
but adding edges $\lineseg {v_{n-1}}{v_0}$ of length 1 
and $\lineseg {v_{n-1}}{v_i}$ of length $\sqrt 2$ for $1\le i\le n-2$,
if these are not already present.
If $\varphi \in \conf(\ul)$, then $\varphi (v_{n-1})$ lies on the circle
of radius 1 about 0 in the $x_{n-1}x_n$ plane.
So there is a unique $g\in SO(2)$ so that $g^{-1}\varphi \in \conf(\ul)$,
and $g$ is a polynomial function of $\varphi (v)$.
So we get an isomorphism 
$\beta ' \colon \conf(\ul')\times SO(2)\to \conf(\ul)$
given by $\beta ' (\varphi ,g )=g \varphi $.
Define $\alpha \colon Z\times F\to ({\mathbb  R}^n)^m$  by 
$\alpha (x,f)=(\gamma (x)e_{n-1},\gamma (x)\alpha '(\gamma (x)^{-1}x,f))$.
Note that the map $(x,f)\mapsto (x,\alpha (x,f)$ is a composition of the 
analytic isomorphisms $\beta ^{-1}\times id\colon Z\times F\to Z'\times SO(2)\times F$ and
$(z,g,f)\mapsto ((z,e_{n-1},\alpha '(z,f)),g)\in \conf(\ul')\times SO(2)$ and
$\beta '$.
Thus it is an analytic isomorphism from $Z\times F$ to $\conf(\ul)$.

So it only remains to prove the indicated converse of part 2.
So suppose $\ul$ is a connected planar classical linkage with one fixed vertex
$v_1$ fixed at $z_1$, and some other nonfixed vertices
$v_2,\ldots ,v_k$.
After reordering the vertices, we may suppose that
 $\lineseg{v_1}{v_k}$
is an edge of $\ul$.
Let this edge have length $r$.
Let $\ul'$ be obtained from $\ul$ by fixing
$v_k$ at some point $z_k$ with $|z_k-z_1|=r$.
Let $Z'=\conf(\ul')$.
Let $G^+\subset \Eu(2)^+$ be the group of rotations about $z_1$.
Then $Z'\times G^+\to \conf(\ul)$ is an isomorphism since
$\conf(\ul)$ is $G^+$ invariant by Lemma \ref{lem:6.3},
and for any $\varphi \in \conf(\ul)$ there is a unique $g\in G^+$
which rotates $\varphi (v_k)$ to $z_k$,
and this $g$ is a polynomial function of $\varphi (v_k)$.
\end{proof}

\begin{proof}
(of Theorem \ref{thm:5.2})
Let $Z''=\beta (Z'\times SO(2))$.
Note that $Z''$ is $SO(2)$ invariant.
Also, the map $\eta \colon Z''\times \Tran(2)\to Z$ 
given by $\eta (z,\tau )=\tau z$ is an isomorphism,
since any $g\in \Eu(2)^+$ can be uniquely decomposed as
$g=\tau \gamma $ for $\tau \in \Tran(2)$ and $\gamma \in SO(2)$.

By Theorem \ref{thm:5} there is a connected classical linkage $\ul'$
with just one fixed vertex, a finite set $F$, 
and an analytic $\alpha '\colon Z''\times F\to ({\mathbb R}^2)^m$ so that
\begin{equation}\label{eqn:20}
\conf(\ul')=\{(x,\alpha '(x,f))\mid x\in Z'' \mand f\in F\}
\end{equation}
and in fact the map $(x,f)\mapsto (x,\alpha '(x,f))$ is an 
analytic isomorphism from $Z''\times F$ to $\conf(\ul')$.
Moreover, the fixed vertex $v$ is fixed at $0$.
Let $\ul$ be formed from $\ul'$ by unfixing the vertex $v$.
Let $\alpha \colon Z\times F\to ({\mathbb R}^2)^m$ be the
analytic map $\alpha (z,f)=\tau \alpha '(\tau ^{-1}z,f)$ where $\tau \in \Tran(2)$ is the
unique translation so that $\tau ^{-1}z\in Z''$.
To be precise, $\eta ^{-1}(z)=(\tau ^{-1}z,\tau )$.

Note that the map $(z,f)\mapsto (z,\alpha (z,f))$ is  the composition
of analytic isomorphisms $\eta ^{-1}\times id\colon Z\times F\to Z''\times \Tran(2)\times F$
and $(z,\tau ,f)\mapsto ((z,\alpha '(z,f)),\tau )\in \conf(\ul')\times \Tran(2)$
and  $(\varphi ,\tau )\mapsto \tau \varphi \in \conf(\ul)$ (which is an isomorphism by
Lemma \ref{lem:3.4}).
Hence it gives an analytic isomorphism from $Z\times F$ to $\conf(\ul)$
as desired.

So the first part of Theorem \ref{thm:5.2} is proven.
Now suppose that $\ul$ is a connected classical linkage with no fixed vertices
and at least two vertices.
We may then pick two vertices $v$ and $w$ of $\ul$ so that
$\lineseg vw$ is an edge of $\ul$.
Let $\ul'$ be the linkage obtained from $\ul$ by 
fixing $v$ to $0$ and fixing $w$ to some point $z_0$
with $|z_0|=\ell (\lineseg vw)$.
Note that $\conf(\ul')\subset \conf(\ul)$.
Moreover the map $\conf(\ul')\times \Eu(2)^+\to \conf(\ul)$ is an analytic
isomorphism since for any $\varphi \in \conf(\ul)$ we know that
$\gamma \varphi \in \conf(\ul')$ where $\gamma \in \Eu(2)^+$ is the unique
element so $\gamma \varphi (v)=0$ and $\gamma \varphi (w)=z_0$,
and $\gamma $ is a polynomial function of $\varphi (v)$ and $\varphi (w)$.
\end{proof}

\section{Which Functions have Functional Linkages?}

Suppose $\ul$ is a (quasi)functional linkage for some function $f$
with domain $X$.
What functions $f\colon X\to ({\mathbb  R}^n)^m$ are possible?
We can completely characterize \qfl s as follows:

\begin{thm}\label{thm:20}
Suppose $X\subset ({\mathbb  R}^n)^k$ and $f\colon X\to ({\mathbb  R}^n)^m$ is a map.
Then the following are equivalent:
\begin{enumerate}
\item  There is a \qfl\ $\ul$ for $f$ with domain $X$.
\item  There is a classical \qfl\ $\ul$ for $f$ with domain $X$.
\item  The graph of $f$ is a semialgebraic set and 
after perhaps permuting the ${\mathbb  R}^n$ factors we have:
\begin{enumerate}
\item  $X=Y_0\times Y_1\times \cdots Y_\ell $, where $Y_i\subset ({\mathbb  R}^n)^{k_i}$,
(and $k_0=0$ is allowed, but $k_i\ge 1$ for $i\ge 1$).
\item  $f$ is a product of maps $f_i\colon Y_i\to ({\mathbb  R}^n)^{m_i}$,
(where $m_i=0$ is allowed and corresponds to composition with projection).
\item  $Y_0$ is compact.
\item  If $i\ge 1$ then $Y_i$ is invariant under the action of $\Eu(n)$
 with compact quotient.
\item  If $i\ge 1$, $f_i$ is  $\Eu(n)$ equivariant.
  That is, for every $\beta \in \Eu(n)$ we have
    $\beta (f_i(z))=f_i(\beta (z))$.
\end{enumerate}
\end{enumerate}
\end{thm}

\begin{proof}
It is trivial that 2 implies 1.
Let us see why 1 implies 3.
If $q$ and $p$ are the input and output maps of $\ul$, then
the graph of $f$ is the image of the polynomial map
$(q,p)\colon \conf(\ul)\to ({\mathbb  R}^n)^k\times ({\mathbb  R}^n)^m$.
So the graph of $f$ is semialgebraic by \cite{S}.
By Lemma \ref{lem:4.2c} below, we may as well assume that $\ul$ is connected.
If $\ul$ has any fixed vertices, then $\conf(\ul)$ is compact by Theorem \ref{thm:4x},
so $X$ is compact since it is the image of the input map $q$.
So we may then take $\ell =0$ and part 3 will hold true. 
On the other hand, if $\ul$ has no fixed vertices then
$\conf(\ul)$ is invariant under the action of $\Eu(n)$ with compact
quotient by Theorem \ref{thm:4x}.
We take $\ell =1$ and $Y_0=$ a point.
Note that the input and output maps $q$ and $p$ are $\Eu(n)$ equivariant.
Consequently $Y_1=q(\conf(\ul))$ is invariant under the action of 
$\Eu(n)$ with compact
quotient
(since $q$ induces a continuous map from $\conf(\ul)/\Eu(n)$
onto $Y_1/\Eu(n)$).
Note that $f$ is $\Eu(n)$ equivariant since if $\beta \in \Eu(n)$ we have
$$f(\beta z)=f(\beta q(y))=f(q(\beta y))=p(\beta y)=\beta p(y)=\beta f(q(y))=\beta f(z)$$
for all $z\in Y_1$ and $y\in q^{-1}(z)$.
So the implication 1 implies 3 is shown.

Now let us show that 3 implies 2.
Note that the graph $G_i$ of $f_i$ is a semialgebraic set since it is 
a projection of the graph of $f$.
If $i>0$ then by 3d) and 3e), the graph $G_i$ is invariant
under the action of $\Eu(n)$ with compact quotient.
By 3c), $G_0$ is compact.
Consequently, by Theorem \ref{thm:4}
there are  classical linkages $\ul_i$ and $U_i\subset \vrtcs(\ul_i)$
so that $\sconf(\ul_i,U_i)=G_i$.
Since $Y_i$ is a coordinate projection of $G_i$,
there is a  $W_i\subset U_i$ 
so that $\sconf(\ul_i,W_i)=Y_i$.
Letting $W_i$ be the input vertices and $U_i-W_i$ be the output vertices,
we thus get a \qfl\ for $f_i$ with domain $Y_i$.
By 3a) and 3b) and Lemma \ref{lem:6.9} we then see that
the disjoint union of the $\ul_i$ is a \qfl\ for $f$.
\end{proof}

Now that \qfl s are classified, we attempt to classify
functional linkages.
This is a bit trickier.
For example, the function $x\mapsto |x_1|e_1$ with domain the 
cube $[-1,1]^n$ has a \qfl\ since its graph is compact and semialgebraic.
But it has no functional linkage whose
restricted domain includes a neighborhood of the point 0.
This is because if $\ul$ is a functional linkage for $f$ with \rd\ $U$,
then $f|_U$ must be an analytic function, since if 
$q$ and $p$ are the input and output maps, there is a finite set $F$
and an analytic 
isomorphism $\sigma  \colon U\times F\to q^{-1}(U)$ so that
$q(\sigma  (u,c))=u$ for all $c\in F$.
Since $\ul$ is functional, we know that $p(\sigma  (u,c))=f(u)$
and thus $f$ is the composition of two analytic functions.

The following two theorems are restricted to the compact case because we
don't know fine enough information about noncompact
configuration spaces.

\begin{thm}\label{thm:21}
Suppose $X\subset ({\mathbb  R}^n)^k$ is compact 
and $f\colon X\to ({\mathbb  R}^n)^m$ is a map,
and $U\subset X$.
Then the following are equivalent:
\begin{enumerate}
\item  There is a functional classical linkage $\ul$ for $f$ with domain $X$
and restricted domain $U$.
\item  There are a compact real algebraic set $Y$, and polynomial maps
$q\colon Y\to X$ and $p\colon Y\to ({\mathbb  R}^n)^m$, a finite set $F$ and an analytic
map $\sigma  \colon U\times F\to Y$ so that:
\begin{enumerate}
\item   $p=f\compose q$.
\item   $q$ is onto.
\item  $q\sigma  (x,c)=x$ for all $(x,c)\in U\times F$.
\item  $\sigma  $ is an analytic isomorphism onto $q^{-1}(U)$.
\end{enumerate}
\end{enumerate}
\end{thm}

\begin{proof}
The implication 1 implies 2 follows immediately from the definition of functional
linkage and the fact that configuration spaces of classical linkages
are algebraic sets.

Let us now see why 2 implies 1.
So suppose we have $Y$, $p$, $q$, $F$, and $\sigma  $ as above.
By replacing $Y$ with the graph of $(p,q)$ we may as well assume that
$p$ and $q$ are given by coordinate projection.
So $Y\subset ({\mathbb R}^n)^m\times ({\mathbb R}^n)^k\times ({\mathbb R}^n)^b$
and $p$ and $q$ are induced by projections to the first and second
batches of coordinates.
By Theorem \ref{thm:1}, there is a classical linkage $\ul$,
a finite set $G$ and an analytic $\beta \colon Y\times G\to ({\mathbb R}^n)^\ell $
so that
$$\conf(\ul)=\{ (y,\beta (y,c)) \mid y\in Y \mand c\in G\}$$
and so the map $(y,c)\mapsto (y,\beta (y,c))$ is an isomorphism
from $Y\times G$ to $\conf(\ul)$.
Then $\ul$ is functional for $f$ with \rd\ $U$.
The first $m$ vertices are the output vertices,
the next $k$ are the input vertices,
and the map $\gamma \colon U\times F\times G\to \conf(\ul)$ giving the trivial analytic cover is
$\gamma (x,c,d)=(\sigma  (x,c),\beta (\sigma  (x,c),d))$.
\end{proof}

\begin{thm}\label{thm:22}
Suppose $X\subset ({\mathbb  R}^n)^k$ is compact 
and $f\colon X\to ({\mathbb  R}^n)^m$ is a map.
Then the following are equivalent:
\begin{enumerate}
\item  There is a strong functional linkage $\ul$ for $f$ with domain $X$.
\item  There is a \qas\ $Y$, polynomial maps
$q\colon Y\to X$ and $p\colon Y\to ({\mathbb  R}^n)^m$ so that:
\begin{enumerate}
\item  $p=f\compose q$.
\item  $q$ is an analytic isomorphism onto $X$.
\end{enumerate}
\end{enumerate}
\end{thm}

\begin{proof}
One direction follows immediately from the definition of strong functional
linkage and the fact that configuration spaces of \rl s
are quasialgebraic sets.

So suppose we have $Y$, $p$, and $q$ as above.
By replacing $Y$ with the graph of $(p,q)$ we may as well assume that
$p$ and $q$ are given by coordinate projection.
So  $Y\subset ({\mathbb R}^n)^m\times ({\mathbb R}^n)^k\times ({\mathbb R}^n)^b$
and $p$ and $q$ are induced by projections to the first and second
batches of coordinates.
By Theorem \ref{thm:2}, there is a \rl\ $\ul$
 and an analytic $\beta \colon Y\to ({\mathbb R}^n)^\ell $
so that
$$\conf(\ul)=\{ (y,\beta (y)) \mid y\in Y\}$$
Then $\ul$ is strongly functional for $f$ with domain $X$.
The first $m$ vertices are the output vertices,
the next $k$ are the input vertices.
The input map is projection to $Y$, followed by $q$
and is hence an analytic isomorphism.
\end{proof}

Theorems \ref{thm:21} and \ref{thm:22} were restricted to the compact case
and so did not completely classify functional linkages.
In the planar case, however, we have enough information
to completely classify functional classical linkages.

\begin{thm}\label{thm:25}
Suppose $X\subset ({\mathbb  R}^2)^k$ 
and $f\colon X\to ({\mathbb  R}^2)^m$ is a map,
and $U\subset X$is nonempty and open.
Then the following are equivalent:
\begin{enumerate}
\item  There is a functional classical linkage 
$\ul$ for $f$ with domain $X$
and restricted domain $U$.
\item  After perhaps permuting the ${\mathbb  R}^2$ factors we have:
\begin{enumerate}
\item  $X=X_0\times X_1\times \cdots \times X_\ell $, 
where $X_i\subset ({\mathbb  R}^2)^{k_i}$,
(and $k_0=0$ is allowed, but $k_i\ge 1$ for $i\ge 1$).
\item  $f$ is a product of maps $f_i\colon X_i\to ({\mathbb  R}^2)^{m_i}$,
(where $m_i=0$ is allowed and corresponds to composition with projection).
\item For $i=0,\ldots ,\ell$ there are real algebraic sets $Y_i$ and
polynomial maps $q_i\colon Y_i\to X_i$
and $p_i\colon Y_i\to ({\mathbb R}^2)^{m_i}$  and a finite set $F$
and an analytic isomorphism $\sigma \colon U\times F\to (q_0\times \cdots \times q_\ell )^{-1}(U)$
so that:
\begin{enumerate}
\item $p_i=f_i\compose q_i$
\item $q_i$ is onto.
\item $(q_0\times \cdots \times q_\ell )(\sigma (x,c))=x$ for all $(x,c)\in U\times F$.
\end{enumerate}
\item  $Y_0$ is compact.
\item  If $i\ge 1$ then $X_i$ and $Y_i$ are invariant under the action of $\Eu(2)$
 with compact quotient.
 \item If $i\ge 1$ then $p_i$ and $q_i$ are $\Eu(2)$ equivariant.
\item If
$i\ge 1$ and $m_i>0$ then there is a compact real algebraic subset $Y_i'\subset Y_i$
so that the map $(y,g)\mapsto gy$ is an isomorphism from
$Y_i'\times \Eu(2)^+$ to $Y_i$.
\end{enumerate}
\end{enumerate}
\end{thm}

\begin{proof}
First let us prove that 2 implies 1.
By replacing each $Y_i$ by the graph of $(q_i,p_i)$, we may as well assume that
$p_i$ and $q_i$ are given by coordinate projections.
By  Theorems \ref{thm:5} and  \ref{thm:5.2}, there are connected classical linkages
$\ul_i$, finite sets $F_i$, and analytic 
$\alpha _i\colon Y_i\times F_i\to ({\mathbb R}^2)^{\ell _i}$ so that
\begin{equation}\label{eqn:5}
\conf(\ul_i)=\{(x,\alpha _i(x,c)) \mid x\in Y_i \mand c\in F_i\}
\end{equation}
and the map $(x,c)\mapsto (x,\alpha _i(x,c))$ is an analytic isomorphism.
Moreover, if $i\ge 1$,
 $\ul_i$ has no fixed vertex.
Let $\ul$ be the disjoint union of the $\ul_i$.
Recall by Lemma \ref{lem:6.2} that $\conf(\ul)=\conf(\ul_0)\times \cdots \times \conf(\ul_\ell )$.

The  projection maps $\pi _i\colon \conf(\ul_i)\to Y_i$ are trivial analytic coverings.
Hence if $\pi =\pi _0\times \cdots \times \pi _\ell $, then $\pi \colon \conf(\ul)\to Y_0\times \cdots \times Y_\ell $
is a trivial analytic covering.
So if $q\colon \conf(\ul)\to X$ is defined by $q=(q_0\times \cdots \times q_\ell )\compose \pi $,
we know by c) that $q$ restricts to a trivial analytic covering
$q|\colon q^{-1}(U)\to U$.

Let $p=(p_0\times \cdots \times p_\ell )\compose \pi $.
Since $p_i$ and $q_i$ are given by coordinate projections from $Y_i$,
we know $p_i\compose \pi $ and $q_i\compose \pi $ are given by coordinate projections from
$\conf(\ul)$.
Consequently,  there are vertices $w_1,\ldots ,w_k$
and $v_1,\ldots ,v_m$ of $\ul$
so that
\begin{eqnarray*}
q &=& (\rh{}{w_1},\ldots ,\rh{}{w_k})\\
p &=& (\rh{}{v_1},\ldots ,\rh{}{v_m})
\end{eqnarray*}
Let the $w_i$ be input vertices and the $v_i$ be output vertices.
Then the input and output maps are $q$ and $p$ respectively.
Note that 
$$p =p_0\compose \pi _0\times \cdots \times p_\ell \compose \pi _\ell 
=f_0\compose q_0\compose \pi _0\times \cdots \times f_\ell \compose q_\ell \compose \pi _\ell =f\compose q$$
 so $\ul$ is quasifunctional for $f$.
Also the domain of $\ul$ is $q(\conf(\ul))=q_0(Y_0)\times \cdots \times q_\ell (Y_\ell )=X$ by 2(c)(ii).
Consequently $\ul$ is a functional linkage with domain $X$ and restricted domain $U$.
So the implication 2 implies 1 is shown.

Now let us show that 1 implies 2.
Let $\ul$ be a functional classical linkage for $f$
with domain $X$ and restricted domain $U$.
We may write $\ul$ as the disjoint union of
$\ul_i$ for $i=0,\ldots ,\ell$ so that
 $\ul_0$ is the union of components of $\ul$ which have fixed vertices
 and so each $\ul_i$ for $i\ge 1$ is connected and has no fixed vertices.
 Let $Y_i=\conf(\ul_i)$.
  
 By Lemma \ref{lem:4.2c} below, we know that each 
 $\ul_i$ is quasifunctional for some $f_i\colon X_i\to ({\mathbb R}^2)^{m_i}$
 and that 2a) and b) hold.
Let $q$ and $p$ be the input and output maps for $f$.
Note that $q=q_0\times \cdots \times q_\ell $ and $p=p_0\times \cdots \times p_\ell $.
Then $p=f\compose q$ implies that $p_i=f_i\compose q_i$ for each $i$.
Also each $q_i$ maps onto $X_i$ since $q$ maps onto $X$.
 Note that 2(c)(iii) holds by definition of \rd.
 By Theorem \ref{thm:4x} we know $X_0$ is compact since it is
 the projection of the compact $\conf(\ul_0)$.
 By part 1 of Theorem \ref{thm:4x} we know 2e),
 and 2f) follows as in the proof of Theorem \ref{thm:20}.
 Finally  2g) follows from Theorem \ref{thm:5.2}.
\end{proof}

We used above the following converse to Lemma \ref{lem:6.9}.

\begin{lem}\label{lem:4.2c}
Suppose $\ul$ is a \qfl\ for $f\colon ({\mathbb R}^n)^k\to ({\mathbb R}^n)^m$ and $\ul$
is the disjoint union of two linkages $\ul_0$ and $\ul_1$.
Then each $\ul_i$ is a \qfl\ for some 
$f_i\colon ({\mathbb R}^n)^{k_i}\to ({\mathbb R}^n)^{m_i}$
and after perhaps reordering the coordinates, $f=f_0\times f_1$.
The input and output maps for $f$ are also the products of the
input and output maps for the $f_i$.

Moreover,  if $\ul$ is a  (strong) functional linkage with \rd\ $U$
and $U$ is a nonempty open subset of the domain of $f$,
then the $\ul_i$ will also be (strong) functional linkages
with \rd\ $U_i$, with $U_i$ nonempty and open in the domain of $f_i$.
\end{lem}

\begin{proof}
Let $w_1,\ldots ,w_k$ and $v_1,\ldots ,v_m$ be the input and output
vertices of $\ul$.
Reorder these so that $w_i\in \ul_0$ if and only if $i\le k_0$ and
$v_i\in \ul_0$ if and only if $i\le m_0$.
Let $k_1=k-k_0$ and $m_1=m-m_0$.
Recall that $\conf(\ul)=\conf(\ul_0)\times \conf(\ul_1)$ by Lemma \ref{lem:6.2}.
Let $q\colon \conf(\ul_0)\times \conf(\ul_1)\to ({\mathbb R}^n)^{k_0}\times ({\mathbb R}^n)^{k_1}$ be the input map.
Note that $q=q_0\times q_1$ for some $q_i\colon \conf(\ul_i)\to ({\mathbb R}^n)^{k_i}$
and likewise the output map $p$ of $\ul$ is $p_0\times p_1$.
Choose any $\varphi '_i\in \conf(\ul_i)$ and let $q(\varphi '_0,\varphi '_1)=(x_0',x_1')$.
If $\ul$ is a (strong) functional linkage with nonempty open \rd\ $U$,
choose them so that $(x'_0,x'_1)\in U$.

We may write $f$ as $f(x_0,x_1)=(f'_0(x_0,x_1),f'_1(x_0,x_1))$
where $x_i\in ({\mathbb R}^n)^{k_i}$ and $f'_i(x_0,x_1)\in ({\mathbb R}^n)^{m_i}$.
Since $f\compose (q_0\times q_1)=p_0\times p_1$ we know that $f'_i\compose (q_0\times q_1)=p_i$.
So for any $\varphi _i\in \conf(\ul_i)$,
$$f'_0(q_0(\varphi _0),q_1(\varphi _1))=p_0(\varphi _0)=f'_0(q_0(\varphi _0),x_1')$$
and thus $f'_0(x_0,x_1)$ is independent of the $x_1$ coordinate
and we may define $f_0(x_0)=f_0'(x_0,x_1')$.
Similarly, $f'_1$ is independent of $x_0$
and we may define $f_1(x_1)=f_1'(x_0',x_1)$.
Then $f_i\compose q_i=p_i$ and $\ul_i$ is a \qfl\ for $f_i$.
Also $f=f_0\times f_1$.

Note that the domain of $f$ is the cartesian product of the domains
of $f_i$.
So if $\ul$ is a (strong) functional linkage with open \rd\ $U$, then
after restricting the domain further, we may as well suppose
that $U=U_0\times U_1$ for $U_i\subset ({\mathbb R}^n)^{k_i}$,
and $x_i'\in U_i$ and $U_i$ open in the domain of $f_i$.
In the strong case, we also may as well assume that $U$ and $U_i$
are the domains of $f$ and $f_i$.
In the nonstrong case, we may restrict each $U_i$ further and assume it is 
connected.
We may do this because the domain of $f_i$ is a semialgebraic set,
hence it is locally
connected by, for example, the curve selection lemma
\cite{M} or triangulability \cite{L}.

By (strong) functionality, there is a finite set $F$ and an analytic isomorphism
$$\sigma  \colon U_0\times U_1\times F\to q^{-1}(U_0\times U_1)
=q_0^{-1}(U_0)\times q_1^{-1}(U_1)$$
so that $q\sigma  (x_0,x_1,c)=(x_0,x_1)$ and so that $F$ is a single point
in the strong case.

Let $F_0\subset F$ be the subset so that 
$x_0'\times x_1'\times F_0=\sigma  ^{-1}(q_0^{-1}(x_0')\times \varphi _1')$.
By connectedness of $U_0$ and discreteness of $F$
(in the nonstrong case), or by the fact that $F_0=F=$ a point
(in the strong case), we then know that
$x_0\times x_1'\times F_0=\sigma  ^{-1}(q_0^{-1}(x_0)\times \varphi _1')$ for any $x_0\in U_0$.
We thus get an analytic isomorphism 
$\sigma  _0\colon U_0\times F_0\to q_0^{-1}(U_0)$ given by
$(\sigma  _0(u,c),\varphi _1')=\sigma  (u,x_1',c)$.
So $\ul_0$ is a (strong) functional linkage for $f_0$ with \rd\ $U_0$.
Similarly, $f_1$ is a (strong) functional linkage with \rd\ $U_1$.
%
%
%
%
\end{proof}

\begin{rem}
Note that in Lemma \ref{lem:4.2c} we have $U_0\times U_1\subset U$ and in fact,
we may choose $U_i$ so that $U_0\times U_1$ contains any one given point of $U$.
\end{rem}

\section{Simulating Higher Dimensional Linkages}

Suppose $K$ is a finite simplicial complex and we identify each simplex of
$K$ with a particular (linear) simplex in some euclidean space.
We could then study the configuration space of all realizations
of $K$ in some ${\mathbb R}^n$, i.e., all maps of $K$ to ${\mathbb R}^n$ 
which restrict to
an isometry on each simplex of $K$.
However,  by Lemma \ref{lem:e} and Lemma \ref{lem:f} below,
 a simplex in ${\mathbb R}^n$ is determined up to
Euclidean motions by the lengths of its edges.
Thus the configuration space of $K$ is the same as the
configuration space of its one dimensional skeleton,
the union of its vertices and edges.
Consequently by looking at higher dimensional
simplicial linkages, we get nothing different
in the way of configuration spaces.

We could generalize even further and still get no new
configuration spaces.
We could consider configuration spaces of polyhedra glued together,
but we could still simulate such objects by a one dimensional linkage.
In particular, let $K_i$, $i=1,\ldots m$ be realizations of finite simplicial
complexes in some Euclidean space.
Let $K_{ji}\subset K_i$ be possibly empty subcomplexes, $j=1,\ldots ,i-1,i+1,\ldots ,m$,
and let $\varphi _{ji}\colon K_{ji}\to K_{ij}$ be compatible simplicial isometries.
In particular,  $\varphi _{ji}=\varphi _{ij}^{-1}$,
$\varphi _{ji}(K_{ji}\cap K_{ki})=K_{ij}\cap K_{kj}$ and
$\varphi _{kj}(\varphi _{ji}(x))=\varphi _{ki}(x)$ for all $x\in K_{ij}\cap K_{kj}$.
Since $\varphi _{ji}$ is simplicial, it takes vertices to vertices, edges to edges, etc.
We may form an object $K$ by gluing the $K_i$ together
using the maps $\varphi _{ji}$.
We may then look at the configuration space of maps from $K$
to ${\mathbb R}^n$ which restrict to an isometry on each $K_i$.
For example $K$ could be the surface of a unit cube and the $K_i$
are its faces.

But the configuration space of $K$ is the same as the configuration space
of some one dimensional linkage $K'$.
To construct $K'$, we just replace each $K_i$ with the one dimensional
complex $K_i'$ which has the same vertices as $K_i$ and has one edge 
of the appropriate length for each pair of vertices of $K_i$.
Glue these one dimensional complexes together using the maps
$\varphi _{ji}$ to identify vertices.
Also identify edges between pairs of identified vertices,
 and thus obtain
a one dimensional simplicial complex $K'$.
In our unit cube example we would get a one dimensional 
complex consisting of the edges
of the unit cube  and the two diagonals on each of its faces.
By Lemma \ref{lem:f} below, the isometries of $K_i$ and $K_i'$ are the same,
so the resulting configuration spaces of $K$ and $K'$ are the same.

Of course if we generalize further to higher dimensional linkages
with curved faces, such reductions are no longer possible.

Following are Lemmas proving various results referred to above.
No doubt these are well known.

\begin{lem}\label{lem:d}
Suppose $z_i\in {\mathbb R}^n$, $i=1,\ldots k$.
Let $T$ be the affine span of $z_1,\ldots ,z_k$.
Then for any fixed $z\in T$ the distance from $z$ to any point $x\in {\mathbb R}^n$
is determined by the distances from the $z_i$ to $x$.
\end{lem}

\begin{proof}
We may write $z=\sum _{i=1}^k t_iz_i$ where $\sum _{i=1}^k t_i =1$.
Then
\begin{eqnarray*}
|x-z|^2 & = &| \sum _{i=1}^k t_i(x-z_i) |^2\\
&=&  \sum _{i,j} t_it_j (x-z_i)\cdot (x-z_j)
\end{eqnarray*}
But $(x-z_i)\cdot (x-z_j)=(|x-z_i|^2+|x-z_j|^2-|z_i-z_j|^2)/2$.
\end{proof}

\begin{lem}\label{lem:g}
Suppose $K\subset {\mathbb R}^m$ and 
suppose $\alpha \colon K\to {\mathbb R}^n$ is an affine transformation
which restricts to an isometry on $K$.
Then $\alpha $ restricts to an isometry on the affine span of $K$.
\end{lem}

\begin{proof}
If $S$ is the affine span of $K$ and $x,y\in S$ then there are
$x_i\in K$, $i=1,\ldots ,k$ and $t_i$ and $s_i$ so that
$x=\sum _{i=1}^k t_ix_i$, $y=\sum _{i=1}^k s_ix_i$, and
$\sum _{i=1}^k t_i=\sum _{i=1}^k s_i=1$.
Since $\alpha $ restricts to an isometry on $K$, we know that
$(\alpha x_i-\alpha x_\ell )\cdot (\alpha x_j-\alpha x_\ell )=(x_i-x_\ell )\cdot (x_j-x_\ell )$
for all $i,j,\ell $.   So
\begin{eqnarray*}
|\alpha x-\alpha y|^2  & = &  |\sum _{i=1}^k (t_i-s_i)(\alpha x_i-\alpha x_1)|^2\\
   & = &  \sum _{i,j} (t_i-s_i)(t_j-s_j) (\alpha x_i-\alpha x_1)\cdot (\alpha x_j-\alpha x_1)\\
    & = &   |\sum (t_i-s_i)(x_i-x_1)|^2 =  |x-y|^2
\end{eqnarray*}
\end{proof}

\begin{lem}\label{lem:e}
Let $K\subset {\mathbb R}^m$ be the realization of some compact simplicial complex
and suppose $\varphi _i\colon K\to {\mathbb R}^n$, $i=0,1$ are two isometric embeddings.
Then there is a $\beta \in \Eu(n)$ so that $\varphi _1=\beta \compose \varphi _0$.
There is also an affine transformation $\alpha \colon {\mathbb R}^m\to {\mathbb R}^n$ so that
$\varphi _0$ is the restriction of $\alpha $.
\end{lem}

\begin{proof}
Suppose first that $K$ is zero dimensional,
so $K$ is a finite number of points $x_i\in {\mathbb R}^m$, $i=1,\ldots ,k$.
By induction on $k$, after composing $\varphi _1$ with some $\beta '\in \Eu(n)$
we may as well assume that 
there is an affine transformation $\alpha '\colon {\mathbb R}^m\to {\mathbb R}^n$
so that $\alpha '(x_i)=\varphi _0(x_i)=\varphi _1(x_i)$
for all $i\ne k$.
Let $T$ be the affine span of $\varphi _0(x_1),\ldots ,\varphi _0(x_{k-1})$.
By Lemma \ref{lem:d} we know that
$\varphi _0(x_k)$ and $\varphi _1(x_k)$ have the same distance to each point
of $T$.
In particular, their orthogonal projections to $T$ coincide since
the orthogonal projection is the point on $T$ of minimum distance.
So we may find a $\beta \in \Eu(n)$ which fixes $T$ and takes
$\varphi _0(x_k)$ to $\varphi _1(x_k)$.
(If $z_0$ is their common orthogonal projection,
choose a rotation $Q\in O(n)$ so $Q(\varphi _0(x_k)-z_0)=\varphi _1(x_k)-z_0$
and so $Q$ fixes $(\varphi _0(x_k)-z_0)^\perp \cap (\varphi _1(x_k)-z_0)^\perp $.
Let $\beta (z)=z_0+Q(z-z_0)$.)

Let us now construct $\alpha $.
Since $\alpha '$ restricts to an isometry on $x_1,\ldots x_{k-1}$,
 it restricts to an isometry of the affine span $S$ of $x_1,\ldots x_{k-1}$
 by Lemma \ref{lem:g}.
So if $x_k\in S$ then  $\alpha '$ restricts to an isometry of $K$ and hence 
$\alpha '(x_k)=\varphi _0(x_k)$ by Lemma \ref{lem:d}.
If $x_k\not\in S$,
we may by hand pick an affine $\alpha $ so $\alpha $ restricts to
$\alpha '$ on $S$ but $\alpha (x_k)=\varphi _0(x_k)$.
In particular, if $\alpha (t(x_k-x)+y)=t(\varphi _0(x_k)-\alpha '(x))+\alpha '(y)$
for all $x,y\in S$ and $t\in {\mathbb R}$, this defines $\alpha $ on the affine span of
$x_1,\ldots ,x_k$ and you just extend $\alpha $ in any fashion to an affine map
on all of ${\mathbb R}^n$.

Now suppose $K$ has dimension $k>0$.
By induction on $k$, we may as well assume that 
there is an affine transformation $\alpha \colon {\mathbb R}^m\to {\mathbb R}^n$
so that $\alpha (x)=\varphi _0(x)=\varphi _1(x)$
for all $x$ in the $k-1$ skeleton of $K$.
Take any point $x$ in the interior of a
$k$ simplex $\sigma $ of $K$, and let $v_0,\ldots ,v_k$ be the
vertices of $\sigma $,
so $x=\sum _{i=0}^k t_iv_i $ and $\sum _{i=0}^k t_i=1$.
Let $x'=(x-t_kv_k)/(1-t_k)$, so $x$ is on the line segment
from $x'$ to $v_k$.
Then since $\varphi _i$ is an isometry we know that
$$\varphi _i(x)=t_k\varphi _i(v_k)+(1-t_k)\varphi _i(x')=t_k\alpha (v_k)+(1-t_k)\alpha (x')=\alpha (x)$$
\end{proof}

\begin{lem}\label{lem:f}
Let $K\subset {\mathbb R}^m$ be the realization of some compact simplicial complex
and let $K'\subset K$ be the set of vertices of $K$.
Then any isometry $\varphi \colon K'\to {\mathbb R}^n$ is the restriction of a unique
isometry $\varphi '\colon K\to {\mathbb R}^n$.
\end{lem}

\begin{proof}
By Lemma \ref{lem:e} there is an affine $\alpha \colon {\mathbb R}^m\to {\mathbb R}^n$ which restricts to
$\varphi $ on $K'$.
But $K$ is contained in the affine span of $K'$, so
by Lemma \ref{lem:g}, $\alpha $ restricts to an isometry of $K$.
Uniqueness follows from Lemma \ref{lem:e}.
\end{proof}

\bibliographystyle{amsplain}

\end{document}